\documentclass{article}

\usepackage{enumerate}
\usepackage{amsmath,xspace,amssymb,mathrsfs}

\input xy
\xyoption{all}
\xyoption{2cell}
\UseAllTwocells
\CompileMatrices

\title{A Quillen model structure for Gray-categories}
\author{Stephen Lack%
\thanks{The support of the Australian Research Council and
DETYA is gratefully acknowledged.}
\\School of Computing and Mathematics\\
University of Western Sydney\\
Locked Bag 1797 Penrith South DC NSW 1797\\
Australia\\
email: {\tt s.lack@uws.edu.au}}
\date{}

\newcommand{\A}{{\ensuremath{\mathscr A}}\xspace}
\newcommand{\B}{{\ensuremath{\mathscr B}}\xspace}
\newcommand{\C}{{\ensuremath{\mathscr C}}\xspace}
\newcommand{\D}{{\ensuremath{\mathscr D}}\xspace}
\newcommand{\E}{{\ensuremath{\mathscr E}}\xspace}

\newcommand{\K}{{\ensuremath{\mathscr K}}\xspace}

\newcommand{\V}{{\ensuremath{\mathscr V}}\xspace}
\newcommand{\hoV}{\ensuremath{\textnormal{ho}\V}\xspace}
\newcommand{\W}{{\ensuremath{\mathscr W}}\xspace}
\newcommand{\Gray}{\textnormal{\bf Gray}\xspace}
\newcommand{\hoGray}{\textnormal{ho}\textbf{Gray}\xspace}
\newcommand{\tricat}{\textnormal{\bf Tricat}\xspace}
\newcommand{\hoGrayCat}{\textnormal{\bf \ho\Gray-Cat}\xspace}
\newcommand{\VCat}{\textnormal{\bf \V-Cat}\xspace}
\newcommand{\WCat}{\textnormal{\bf \W-Cat}\xspace}
\newcommand{\graycat}{\textnormal{\bf Gray-Cat}\xspace}
\newcommand{\graygpd}{\textnormal{\bf Gray-Gpd}\xspace}
\newcommand{\SSet}{\textnormal{\bf SSet}\xspace}
\newcommand{\set}{\textnormal{\bf Set}\xspace}
\newcommand{\Top}{\textnormal{\bf Top}\xspace}
\newcommand{\cat}{\textnormal{\bf Cat}\xspace}
\newcommand{\gpd}{\textnormal{\bf Gpd}\xspace}

\newcommand{\bicat}{\textnormal{\bf Bicat}\xspace}

\newcommand{\twocat}{\textnormal{\bf 2-Cat}\xspace}
\newcommand{\twogpd}{\textnormal{\bf 2-Gpd}\xspace}
\newcommand{\threecat}{\textnormal{\bf 3-Cat}\xspace}
\newcommand{\cmptd}{\textnormal{\bf Cmptd}\xspace}
\newcommand{\sesquicat}{\textnormal{\bf SesquiCat}\xspace}

\renewcommand{\AA}{{\ensuremath{\mathbb A}}\xspace}
\newcommand{\BB}{{\ensuremath{\mathbb B}}\xspace}
\newcommand{\EE}{{\ensuremath{\mathbb E}}\xspace}
\newcommand{\TT}{{\ensuremath{\mathbb T}}\xspace}
\newcommand{\PB}{\ensuremath{\mathbb P\mathbb B}\xspace}
\newcommand{\ho}{\textnormal{ho}\xspace}
\newcommand{\ob}{\textnormal{ob}\xspace}

\newcommand{\bb}{\boldsymbol}

\newcommand{\ox}{\otimes}

\newtheorem{theorem}{Theorem}[section]
\newtheorem{lemma}[theorem]{Lemma}
\newtheorem{proposition}[theorem]{Proposition}
\newtheorem{corollary}[theorem]{Corollary}
\newtheorem{preremark}[theorem]{Remark}
\newenvironment{remark}{\begin{preremark}\rm}{\end{preremark}}

\newcommand{\proof}{\noindent{\sc Proof:}\xspace}
\def\endproof{~\hfill$\Box$\vskip 10pt}

\newcommand{\two}{\ensuremath{\mathbf{2}\xspace}}

\begin{document}

\label{firstpage}
\maketitle

\begin{abstract}
A Quillen model structure on the category \graycat of \Gray-categories
is described, for which the weak equivalences are the triequivalences. It is shown to restrict to the full subcategory \graygpd
of \Gray-groupoids. This is used to provide a functorial and model-theoretic 
proof of the unpublished theorem of Joyal and Tierney that
\Gray-groupoids model homotopy 3-types. The model structure on
\graycat is conjectured to be Quillen equivalent to a model structure
on the category \tricat of tricategories and strict homomorphisms of tricategories.
%Construction of model structure on \graycat typed on 3 Feb 2009 based on working note from 8 Jan 2007. Other sections added bit by bit. Very preliminary version: not for further distribution at this stage. 
\end{abstract}

Much of the recent development of higher category theory has been 
inspired and guided by the idea that higher groupoids should classify 
homotopy types. This is sometimes called the ``homotopy hypothesis''.
A higher groupoid is a higher category in which morphisms at all
dimensions are invertible in a suitable sense; here, that will mean that they 
are strictly invertible, so that the composite of a morphism with its inverse 
is literally equal to an identity morphism.
The goal of this paper is to describe certain precise relationships between
higher categories and homotopy types, in low dimensions, using the
machinery of model
categories and Quillen equivalences. 

Experience so far suggests that for each dimension $n$ there should be 
a model category which captures what it is to ``do'' $n$-dimensional 
category theory, and that this should restrict to a model category on 
$n$-dimensional groupoids which provide a model for homotopy $n$-types.

The case $n=1$ is well-understood. The category \cat of (small) categories 
and functors has a model structure \cite{Joyal-Tierney-CatE} for which the weak equivalences are the
equivalences of categories, and the fibrations are the functors with the isomorphism-lifting property. (The word fibration will always be used in the 
model-categorical sense; we shall not need to discuss categorical fibrations.)
This model structure on \cat is sometimes given the 
epithet ``folklore'', but I prefer to call it the {\em categorical model structure},
since it corresponds to what it is to do category theory. In category theory 
one regards two categories as being the same when they are equivalent (and 
then uses the equivalence to identify them).  Furthermore the fibrations for
the model structure arise naturally as the morphisms with the property that
pulling back along them sends equivalences to equivalences; there is also 
an analogous characterizations of the cofibrations \cite{JS-pseudopullbacks}.
This model structure on \cat restricts to the full subcategory \gpd of \cat
consisting of the groupoids. The nerve functor $N:\gpd\to\SSet$
is the right adjoint part of a Quillen functor, and this allows us to regard 
groupoids as models for homotopy 1-types; this is made more precise 
in Section~\ref{sect:localization}, where we see that $N:\gpd\to\SSet$ is 
actually a Quillen equivalence when we replace the usual model structure
on \SSet with a localized version, which kills all homotopy information in
dimension greater than 1. Thus groupoids are seen as the ``intersection 
of category theory and homotopy theory''.
\footnote{We shall not have cause to consider the other well-known model structure
on \cat, due to Thomason \cite{Thomason-Cat} --- this is in fact Quillen equivalent to that of \SSet, and so using this would give a rather different ``intersection''. But while
this Thomason model structure is important in homotopy theory, it seems to correspond to another use of categories than category theory itself.}

Similarly, when we turn to the case $n=2$, the category \twocat of (strict) 2-categories and (strict) 2-functors has a model structure \cite{qm2cat,qmbicat}, for which the weak equivalences are the biequivalences, and the fibrations are 2-functors allowing
liftings of both 1-cells which are equivalences and 2-cells which are invertible.
Once again, these biequivalences capture exactly the usual notion of ``sameness''
for 2-categories, and the fibrations can be characterized in terms of pulling back
of biequivalences. The model structure restricts to the full subcategory \twogpd
of \twocat consisting of the 2-groupoids: these are the 2-categories in which 
all 1-cells and 2-cells are (strictly) invertible. The resulting model structure on 
\twogpd was introduced in \cite{MoerdijkSvensson}, prior to \cite{qm2cat}. Once again there is a nerve functor
$N:\twogpd\to\SSet$ which is the right adjoint part of a Quillen functor, and is
in fact a Quillen equivalence for a localized model structure on \SSet, obtained
by killing homotopy information in dimension greater than $2$; once
again, this last aspect is described in more detail in 
Section~\ref{sect:localization}.

The category \twocat uses strict notions both for objects (2-categories) and morphisms
(2-functors). One might question whether this is really suitable, given that in 
practice many 2-dimensional structures are not strict. The reason for using 2-functors
is that they are much better behaved than the more general pseudofunctors, and
that the more general ones are in any case encoded via the model structure. This 
is the typical situation: use well-behaved gadgets as models, then use the model
structure to get at the more general notions. As far as the objects go, the standard
notion of weak 2-category is called a bicategory, and the category \bicat of bicategories
and strict homomorphisms of these also has a model structure \cite{qmbicat},
which is Quillen equivalent to that on \twocat. This fact includes the result \cite{MacLane-Pare} that every bicategory is biequivalent to a strict one.

When we come to the case $n=3$, the general notion of weak 3-category is 
called a tricategory, and this time it is not the case that every 
tricategory is suitably equivalent (``triequivalent'') to a strict one, however
every tricategory is triequivalent to an intermediate structure called a 
Gray-category \cite{GPS}. One of the main results proved below is that there is
a model structure on the category \graycat of Gray-categories for which the
weak equivalences are the triequivalences. A Gray-groupoid is a Gray-category 
in which all 1-cells, 2-cells, and 3-cells have strict inverses. We shall also 
see that the model structure on \graycat restricts to give one on the full
subcategory \graygpd of Gray-groupoids. There is a nerve functor 
$N:\graygpd\to\SSet$, defined in  \cite{Berger-3type}, 
which turns out to be the right adjoint part of a Quillen adjunction,
and as in the previous cases this becomes a Quillen equivalence when we 
localize \SSet by killing homotopies; this time those in dimension
greater than 3. 

This provides a model-theoretic and functorial formulation of the unpublished but 
widely advertised result of Joyal and Tierney that Gray-groupoids model homotopy 
3-types. The possibility of such a formulation was suggested by Berger in \cite{Berger-3type}. 
In \cite[Theorem~3.3]{Berger-3type}, Berger considered the nerve functor 
$N:\graygpd\to\SSet$ and its left adjoint $\Pi_3$, and asserted that the model 
structure on \SSet could be transported across this adjunction to give a model structure 
on \graygpd, and then show that the adjunction would induce an equivalence 
at the homotopy level between Gray-groupoids and simplicial 3-types. He did not, 
however, give a proof that the model structure actually transports. Once one knows that 
the structure transports, the fact that there is a Quillen adjunction is immediate, while 
the fact that is a Quillen equivalence, and so induces an equivalence of homotopy 
categories, essentially amounts to \cite[Proposition~3.2]{Berger-3type}. 

I take a different approach, constructing the model structure on \graygpd from that on \graycat, rather than transporting it from \SSet. The existence of the model structure
is then more or less immediate, but one must show that the nerve is part of a Quillen adjunction --- this contrasts with
the approach of \cite{Berger-3type}, where the work to be done is in 
constructing the model structure itself. The arguments given here to show that 
the nerve is part of a Quillen adjunction can in fact be used to complete the 
proof in \cite{Berger-3type} of the existence of the model structure. In either approach
one then must do some work (essentially the same) to get a Quillen equivalence:
the key step is \cite[Proposition~3.2]{Berger-3type}.

It turns out, after the fact, that the two model structures on \graygpd actually 
agree: see Corollary~\ref{cor:transport}.

Unlike the
case $n=2$, we do not yet have a model structure on the category of tricategories,
but there is every reason believe that such a structure exists, and that it will
be Quillen equivalent to our model structure on \graycat: we describe what is
known in Section~\ref{sect:tricategories}.

We now give a brief outline of the
paper. Section~\ref{sect:background} recalls some background material on
enriched categories, on model categories, on 2-categories, and on
Gray-categories. Section~\ref{sect:main} describes the model
structure, while Section~\ref{sect:fibrations} gives an alternative
characterization
of the fibrations, and Section~\ref{sect:proof} proves the remaining model category
axioms. In Section~\ref{sect:groupoids}, we look at the restricted
model structure for Gray-groupoids, and its behaviour with respect to
the nerve functor for Gray-groupoids. Section~\ref{sect:localization}
shows how this relationship can be described using Bousfield
localizations of the model category of simplicial
sets. Section~\ref{sect:tricategories} describes briefly what is known
about tricategories, and conjectures what else might be true. The last
two sections concern the cofibrations for the model structure on
\graycat and in particular the cofibrant objects:
Section~\ref{sect:computads} contains some preliminary material on
computads, and Section~\ref{sect:cofibrations} a characterization of 
cofibrations of Gray-categories and of cofibrant Gray-categories.

I first proved the existence of the model structure on \graycat in
January 2007. Since then, the paper has evolved in various ways,
thanks in part to helpful conversations with a number of people.
I received many helpful suggestions from Clemens Berger; among other 
things, when I explained to Clemens the
model structure, he suggested that it might be possible to characterize
the fibrations in the way given below, in terms of $P_*$ and $\pi_*$.
John Harper suggested the particular form that the localizations in 
Section~\ref{sect:localization} might take. I am grateful to both of them, and also to Michael Batanin and Richard Garner, for enlightening discussions on other aspects 
of the paper.

\section{Background}
\label{sect:background}

\subsection{Enriched categories}

For a monoidal category $\V=(\V,\ox,I)$, we write \VCat for the category 
of \V-categories and \V-functors; this can be made into a 2-category with the \V-natural transformations as 2-cells, but we shall not need to do so.
The most important case will be where \V is the category \twocat of 2-categories
and 2-functors, and $\ox$ is what is called the Gray tensor product; then the
\V-categories are precisely the Gray-categories of the introduction. When 
\V is the cartesian closed category \cat, a \V-category is a
2-category. If $A$ and $B$ are objects of a \V-category \A, we write $\A(A,B)$ for
the corresponding hom-object in \V.
%; we shall also use some other monoidal categories \V.

For monoidal categories $\V=(\V,\ox,I)$  and $\W=(\W,\ox,I)$ , a monoidal functor $P:\V\to\W$ consists of a functor $P$ between the underlying categories,
equipped with natural coherent morphisms $PX\ox PY\to P(X\ox Y)$ and 
$I\to PI$. Such a monoidal functor induces a 2-functor
$P_*:\VCat\to\WCat$ 
sending a \V-category \A to the \W-category $P_*\A$ with the same
objects as \A, but with \W-valued homs given by 
$(P_*\A)(A,B)=P(\A(A,B))$. In particular, for any monoidal category \V
the representable functor $\V(I,-):\V\to\set$ arising from the unit $I$
is monoidal, and so induces a 2-functor $\VCat\to\cat$ sending a 
\V-category \A to its underlying ordinary category $\A_0$.

If \V is a monoidal category, a \V-category \A is said to have a 
property ``locally'' if each hom-object $\A(A,B)$ has the property
(as an object of \V). Similarly a \V-functor $F:\A\to\B$ has a property
if each $F:\A(A,B)\to\B(FA,FB)$ has the property (as a morphism of \V).
So for example, one might call a \V-functor $F:\A\to\B$ ``locally an 
isomorphism'' if each $F:\A(A,B)\to\B(FA,FB)$ is an isomorphism in \V;
more commonly, however, such an $F$ is called fully faithful.

A \V-functor $F:\A\to\B$ is an equivalence if there exists a \V-functor
$G:\B\to\A$ with $GF\cong 1$ and $FG\cong 1$. This is the case if and
only if $F$ is fully faithful (in the sense of the previous paragraph)
and essentially surjective on objects (for each $B\in\B$ there is an
$A\in\A$ with $B\cong FA$ in \B). Note that an isomorphism $B\cong FA$
in \B is the same thing as an isomorphism in the underlying
ordinary category $\B_0$ of \B. 

\subsection{Monoidal model categories}

A {\em monoidal model category} \cite{Hovey-book} is a category that 
has both a symmetric monoidal closed structure and a model structure,
satisfying some compatibility conditions: the ``pushout product'' axiom
should hold, as should a condition on the tensor unit; the latter condition
is automatic if, as in our examples, the unit is cofibrant.

All that is really needed in this paper is that 
if \V is a monoidal model category, then the homotopy category \hoV
of \V inherits a derived monoidal structure for which the canonical functor
$P:\V\to\hoV$ is strong monoidal. We shall also write $\pi:\V\to\set$ for
the composite of $P$ with the monoidal functor $\hoV(I,-):\hoV\to\set$ given
by homming out of the unit. This also induces a functor $\pi_*:\VCat\to\cat$. 

\subsection{2-categories}

In this section we recall a few basic ideas about 2-categories. Adjunctions
can be defined in any 2-category, with the usual notion of adjunction being 
the case of the 2-category \cat. An adjunction in a 2-category \K consists
of morphisms $f:A\to B$ and $g:B\to A$ with 2-cells $\eta:1\to gf$ and
$\epsilon:gf\to 1$ satisfying the triangle equations. 

Such an adjunction is called an {\em adjoint equivalence} when the unit $\eta$
and counit $\epsilon$ are both invertible. Occasionally I may allow myself to
say that a morphism $f$ ``is an adjoint equivalence''; this always means that 
some definite choice of $g$, $\eta$, and $\epsilon$ has been made. 

On the other hand, a morphism $f$ is an {\em equivalence} when there {\em exists}
a morphism $g$ with $gf\cong 1$ and $fg\cong 1$. As is well-known
(see \cite{qm2cat} for example), if $f$ is an equivalence and $\eta:1\cong gf$
is an isomorphism, then there is exactly one choice of $\epsilon$ for which the
one (either) triangle equation holds; then the other will also hold, and we shall
have an adjoint equivalence. (Similarly given $\epsilon$ there is exactly one 
way to choose $\eta$.)

The category \twocat is locally finitely presentable, by \cite{vcat}, and so 
small object arguments work smoothly. We shall consider \twocat with
the model structure of \cite{qm2cat}. The corresponding homotopy category
has 2-categories as objects, and pseudonatural equivalence classes of 
pseudofunctors as morphisms. A weak equivalence in \twocat is
(a 2-functor which is) a biequivalence $F:\A\to\B$: this means that $F$ is
locally an equivalence, so that each $F:\A(A,B)\to\B(FA,FB)$ is an equivalence;
and furthermore that $F$ is biessentially surjective, so that each object $X\in\B$
is equivalent to one of the form $FA$ for some $A\in\A$.

\subsection{\Gray-categories}

We shall primarily be interested in the case where \V is the category \twocat,
equipped with the model structure of \cite{qm2cat} and the symmetric monoidal
closed structure given by the Gray tensor product \cite{Gray};  this monoidal category is often called \Gray, to distinguish it from the monoidal category with
the cartesian product as tensor product. 

A 3-category is a category enriched in \twocat with the cartesian product;
a \Gray-category is a category enriched in \twocat with the Gray tensor
product. The difference between \Gray-categories and 3-categories is that 
while given 2-cells 
$$\xymatrix{
A \rtwocell^f_{f'}{\alpha} & B \rtwocell^g_{g'}{\beta} & C}$$
in a 3-category, the composites $\beta f'.g\alpha$
and $g'\alpha.\beta f$ agree, in a \Gray-category they need not; rather,
there is a specified isomorphism which we call simply $\beta\alpha$ or $\beta_\alpha$. %(and if I'm careless I might end up also using $\beta\alpha$ for the inverse).
These are sometimes called the {\em pseudonaturality isomorphisms}.

A morphism $f:A\to B$ in a \Gray-category \AA is called a 
{\em biequivalence in \AA} if there is a morphism $g:B\to A$ with $gf\simeq 1$ 
and $fg\simeq 1$. For any \Gray-functor $F:\AA\to\BB$ we have
$Fg.Ff=F(gf)\simeq F1=1$ and $Ff.Fg\simeq 1$, and so $Ff$ is 
also a biequivalence.

Then \hoGray
is the category of 2-categories and equivalence classes of pseudofunctors,
while $\pi:\Gray\to\set$ sends a 2-category \A to the set of equivalence
classes of objects of \A (in other words objects $A$ and $B$ are equivalent
exactly when they are equivalent in the 2-category \A). We shall write 
\graycat for the category of \Gray-categories and \Gray-functors (the 
2-cells will have little explicit role).
The category \twocat is locally finitely presentable, thus so by  
\cite{vcat} is \graycat. This means that small object arguments work smoothly
in \graycat (all objects are small).

We shall frequently use the functors $P_*:\graycat\to\hoGrayCat$ and
$\pi_*:\graycat\to\cat$.

\section{The model structure on \graycat}
\label{sect:main}

We define a \Gray-functor $F:\AA\to\BB$ to be a 
{\em weak equivalence} if it induces an equivalence 
$P_* F:P_*\AA\to P_*\BB$ of
\hoGray-categories. We shall spell out more explicitly what this means
in the following paragraphs, but it is immediate from this definition that these
weak equivalences are closed under retracts and satisfy the 2-out-of-3
property, since equivalences of enriched categories have these closure
properties.

To say that $P_*F$ is an equivalence is to say that it is fully
faithful and essentially surjective on objects. Being fully 
faithful means that each $P_*F:P_*\AA(A,B)\to P_*\BB(FA,FB)$ is
invertible in \hoGray; in other words, each $PF:P(\AA(A,B))\to P(\BB(FA,FB))$
is invertible in \hoGray; but this is just to say that each
$F:\AA(A,B)\to\BB(FA,FB)$ is a weak equivalence in \Gray;
in other words, a biequivalence.

On the other hand $P_*F$ is essentially surjective when for every 
$B\in\BB$ there is an $A\in\AA$ with an isomorphism $B\cong FA$ in 
$P_*\BB$. But isomorphisms in $P_*\BB$ are just isomorphisms in the 
underlying ordinary category $\pi_*\BB$ of $P_*\BB$, so this is 
just saying that $\pi_*F:\pi_*\AA\to\pi_*\BB$ is essentially surjective
on objects.

We can also make this still more explicit. The category 
$\pi_*\BB$  has the same objects
as \BB, but a morphism in $\pi_*\BB$ from $B$ to $C$ is an equivalence
class of pseudofunctors $1\to\BB(B,C)$. So such a morphism can be represented
by a morphism $f:B\to C$ in \BB, but two such $f,g:B\to C$ represent the
same morphism in $(\ho_*\BB)_0$ if and only if they are equivalent. So 
objects $B$ and $C$ of $\ho_*\BB$ are isomorphic if and only if there exist
$f:B\to C$ and $g:C\to B$ with equivalences $gf\simeq 1$ and $fg\simeq 1$;
that is, if and only if they are {\em biequivalent} in \BB. And this
means that the weak equivalences in \graycat are precisely the {\em triequivalences} of \Gray-categories.

We define a \Gray-functor $F:\AA\to\BB$ to be a {\em fibration} if
it induces fibrations $F:\AA(A,B)\to\BB(FA,FB)$ in \V between the 
hom-objects, for all $A,B\in\AA$, and if moreover the functor  $\pi_*F:\pi_*\AA\to\pi_*\BB$ is a fibration in \cat (an isofibration).
This means that for every object $A\in\AA$ and every isomorphism 
$f:B\cong FA$ in $\pi_*\BB$, there is an object $A'\in\AA$ with $FA'=B$
and an isomorphism $f':A'\cong A$ in $\pi_*\AA$ with $Ff'=f$ (in
$\pi_*\BB$).

Of course a \Gray-functor $F:\AA\to\BB$ is a {\em trivial fibration}
if it is weak equivalence and a fibration; this means that each
$F:\AA(A,B)\to\BB(FA,FB)$ is a weak equivalence and a fibration --- 
that is, a trivial fibration in \twocat, and that each $\pi_* F:\pi_*\A\to\pi_*\B$
is essentially surjective on objects and an isofibration. This certainly
implies that $\pi_*F$ is surjective on objects; but $\pi_*$ does not 
affect the objects, and so $F$ itself is surjective on objects. We can
now prove:

\begin{proposition}
A \Gray-functor $F:\AA\to\BB$ is a trivial fibration if and only if it
is surjective on objects and is locally a trivial fibration in \twocat.
\end{proposition}

\proof
We have already seen the ``only if'' part. For the converse, suppose that
$F$ is surjective on objects and a trivial fibration on the homs. Then
certainly it is a weak equivalence, and is a fibration on the homs. It
remains to show that $\pi_*F$ is an isofibration. But since 
$F$ is a weak equivalence on the homs, $\pi_*F$ is fully faithful; since
it is also surjective on objects, it is a trivial fibration in \cat, 
and so in particular an isofibration.
\endproof

We define the cofibrations to be the \Gray-functors
with the left lifting property with respect to the trivial fibrations.
For a 2-category $X$ we define $\two_X$ to be the \Gray-category with 
two objects $0$ and $1$, with $\two_X(0,0)=\two_X(1,1)=1$, $\two_X(0,1)=X$,
and $\two_X(1,0)=0$. This is clearly functorial in $X$, and for any
\Gray-category \AA, there is a bijection between \Gray-functors
$\two_X\to\AA$ and pairs of objects $A,B\to\AA$ along with a 2-functor
$X\to\AA(A,B)$.  We may now
take as generating cofibrations the $\two_j:\two_X\to\two_Y$ for each
generating cofibration $j:X\to Y$ in \twocat, and the unique
\Gray-functor $0\to 1$; and the cofibrations and trivial fibrations 
from a cofibrantly generated weak factorization system with respect
to these generating cofibrations.

We shall prove the following theorem:

\begin{theorem}\label{thm:main}
This choice of cofibrations, fibrations, and weak equivalences makes
\graycat into a (combinatorial) model category.
\end{theorem}

The key step is to show that fibrations can be defined via a right 
lifting property. We shall do this in the next section. Given that,
the remainder of the proof is largely formal.

\begin{remark}
For any monoidal model category \V we could define a \V-functor 
$F:\AA\to\BB$ to be a weak equivalence if $P_*F:P_*\AA\to P_*\BB$
is an equivalence of \hoV-categories, and to be a fibration if $F$
it is locally a fibration in \V and $\pi_*F:\pi_*\AA\to\pi_*\BB$ is
a fibration in \cat. Then the trivial fibrations would be precisely 
the \V-functors which are surjective on objects and locally a 
trivial fibration in \V, and the cofibrations could be defined to be 
the morphisms with the left lifting property with respect to the trivial
fibrations. In order to prove the model category axioms, however, we
shall use various properties of the monoidal model category $\V=\Gray$.
\end{remark}

% The proof will consist of the following steps:

% \begin{enumerate}[(a)]
% \item Showing that every biequivalence in a \Gray-category can be
% made into an ``adjoint biequivalence''
% \item Showing that fibrations can be defined using an ``adjoint biequivalence
% lifting property'' rather than the ``biequivalence lifting property''
% \item Using this characterization to describe generating trivial 
% cofibrations, and so make the trivial cofibrations and the fibrations
% into a weak factorization system 
% \item Constructing path objects
% \item Using the path objects and the existing factorizations to show
% that the trivial cofibrations are precisely those cofibrations which
% are also weak equivalences
% \end{enumerate}

\section{Fibrations of \Gray-categories and adjoint biequivalences in 
\Gray-categories}
\label{sect:fibrations}

Recall that a \Gray-functor $F:\AA\to\BB$ is a fibration if it satisfies
the following two conditions:

\begin{enumerate}[(i)]
\item The 2-functor $F:\AA(A,B)\to\BB(FA,FB)$ is a fibration in \twocat 
for all objects $A,B\in\AA$
\item For each $A\in\AA$ and each isomorphism $f:B\cong FA$ in $\pi_*\BB$,
there is an $A'\in\AA$ with $FA'=B$ and an isomorphism $f':A'\cong A$ 
in $\pi_*\AA$ with $Ff'=f$.
\end{enumerate}

The first condition clearly amounts to the fact that $F$ has the right 
lifting property with respect to the \Gray-functors $\two_j:\two_X\to\two_Y$
for each generating trivial cofibration $j:X\to Y$ of \twocat. The second 
condition is not, as it stands, a right lifting property of $F$, since 
it involves $\pi_*\AA$ and $\pi_*\BB$. We shall show how to replace 
condition (ii) with a condition (ii*) which is a right lifting property,
and which, in the presence of (i), is equivalent to (ii).

Recall that a morphism $f:A\to B$ in a \Gray-category \BB is a biequivalence
if there exist a $g:B\to A$ with $1\simeq gf$ and $fg\simeq 1$.
We define, following \cite{Verity-thesis}, an {\em adjoint biequivalence} in \BB to consist of $f:A\to B$
and $g:B\to A$, equipped with adjoint equivalences $\eta:1\simeq gf$
and $\epsilon:fg\simeq 1$, and isomorphisms $S:\epsilon f.f\eta\cong 1$
and $T:1\cong g\epsilon.\eta g$ for which the pasting composites
$$\xymatrix{
& gf \ar[dr]_{\eta gf} \ar@/^1pc/[drr]^{1} \drtwocell<\omit>{<-2>Tf} && &
&& fg \ar[dr]^{\epsilon} \\ 
1 \rrtwocell<\omit>{\eta\eta} \ar[ur]^{\eta}\ar[dr]_{\eta} && 
gfgf \ar[r]^{g\epsilon f} & gf &
fg \ar[r]^{f\eta g} \urtwocell<\omit>{<2>fT} \drtwocell<\omit>{<-2>Sg} 
\ar@/^1pc/[urr]^{1} \ar@/_1pc/[drr]_{1} & fgfg \ar[ur]_{fg\epsilon} \ar[dr]^{\epsilon fg} \rrtwocell<\omit>{\epsilon\epsilon} && 1 \\
& gf \ar[ur]^{gf\eta} \ar@/_1pc/[urr]_{1} \urtwocell<\omit>{<2>gS} && &
&& fg \ar[ur]_{\epsilon} }$$
are identities. One might call these last conditions the {\em tetrahedron 
equations}. Clearly $f$ and $g$ are then biequivalences. Note that we
have simply said that $\eta:1\simeq gf$ and $\epsilon:fg\simeq 1$
are adjoint equivalences, but that this must be understood to mean
that we have chosen all the remaining structure of an adjoint
equivalence (e.g. a morphism $\eta^*:gf\to 1$ and invertible 2-cells
$\eta^*\eta\cong1$ and $\eta\eta^*\cong 1$ satisfying the triangle
equations).

We shall show
that every biequivalence can be made into an adjoint biequivalence.
More precisely:

\begin{proposition}\label{prop:adj-biequiv}
Let $f:A\to B$ be a biequivalence. 
\begin{enumerate}[(i)]
\item If $g:B\to A$, and $\epsilon:fg\simeq 1$ is an equivalence,
then there exists an equivalence $\eta:1\simeq gf$ and
an isomorphism $S:\epsilon f.f\eta\cong 1$.
%\item If $g:B\to A$, and $\eta:1\simeq gf$ is an equivalence,
%then there exists an equivalence $\epsilon:fg\simeq 1$ and
%an isomorphism $S:\epsilon f.f\eta\cong 1$.
\item If $g: B\to A$, if $\eta:1\simeq gf$ and $\epsilon:fg\simeq 1$
are equivalences, and if $S:\epsilon f.f\eta\cong 1$ is an
isomorphism, then there exists a unique isomorphism 
$T:g\epsilon.\eta g\cong 1$ satisfying the first tetrahedron equation.
\item If $g:B\to A$, if $\eta:1\simeq gf$ and $\epsilon:fg\simeq 1$
are equivalences, if $S:\epsilon f.f\eta\cong 1$ and 
$T:g\epsilon.\eta g\cong 1$ are isomorphisms, and if the first tetrahedron
equation holds then so does the second.
\end{enumerate}
\end{proposition}

\proof 
% (i) The 2-functor $\AA(f,B):\AA(B,B)\to\AA(A,B)$ is a 
% biequivalence, and so locally an equivalence; in particular the induced functor 
% $$\xymatrix{\relax\AA(B,B)(fg,1) \ar[r]^{\AA(f,B)} & \AA(A,B)(fgf,f) }$$
% is an equivalence and so 
% $$\xymatrix{\relax\AA(B,B)(fg,1) \ar[r]^{\AA(f,B)} & 
% \AA(A,B)(fgf,f) \ar[r]^{\AA(f\eta,f)} & \AA(A,B)(f,f) }$$
% is an equivalence, since $\eta$ is one. So by essential surjectivity,
% there is an $\epsilon:fg\to 1$ and whose image $\epsilon f.f\eta$
% is isomoprhic to the identity on $f$, say by $S:\epsilon f.f\eta\cong 1$.
% It remains to check that $\epsilon$ is an equivalence. Now $\eta$ is 
% an equivalence, so $f\eta$ is an equivalence, so $\epsilon f$ is 
% an equivalence; but $f$ is a biequivalence, so $\epsilon$ is indeed
% an equivalence. 
(i) The 2-functor $\AA(A,f):\AA(A,A)\to\AA(A,B)$ is a 
biequivalence, and so locally an equivalence; in particular the induced functor 
$$\xymatrix{\relax\AA(A,A)(1,gf) \ar[r]^{\AA(A,f)} & \AA(A,B)(f,fgf) }$$
is an equivalence and so 
$$\xymatrix{\relax\AA(A,A)(1,gf) \ar[r]^{\AA(A,f)} & 
\AA(A,B)(f,fgf) \ar[r]^{\AA(f,\epsilon f)} & \AA(A,B)(f,f) }$$
is an equivalence, since $\epsilon$ is one. So by essential surjectivity,
there is an $\eta:1\to gf$ whose image $\epsilon f.f\eta$
is isomorphic to the identity on $f$, say by $S:\epsilon f.f\eta\cong 1$.
It remains to check that $\eta$ is an equivalence. Now $\epsilon$ is 
an equivalence, so $\epsilon f$ is an equivalence, so $f\eta$ is 
an equivalence; but $f$ is a biequivalence, and thus $\eta$ is indeed
an equivalence. 

(ii) The 2-functor $\AA(B,f):\AA(B,A)\to\AA(B,B)$ is a biequivalence,
and so locally an equivalence; in particular the induced functor
$$\xymatrix{\relax
\AA(B,A)(g,g) \ar[r]^-{\AA(B,f)} & \AA(B,B)(fg,fg) }$$
is an equivalence and so 
\begin{equation}\label{e-.f}\xymatrix @C3pc {\relax
\AA(B,A)(g,g) \ar[r]^-{\AA(B,f)} & \AA(B,B)(fg,fg)
\ar[r]^-{\AA(B,B)(fg,\epsilon)} & 
\AA(B,B)(fg,1) }\end{equation}
is an equivalence. Now this composite equivalence maps the identity
on $g$ to $\epsilon$, and maps $g\epsilon.\eta g$ to 
$\epsilon.fg\epsilon.f\eta g$, but we have isomorphisms
\begin{equation}\label{eq1a}\xymatrix{\relax
\epsilon.fg\epsilon.f\eta g \ar[r]^-{\epsilon_\epsilon f\eta g} &
\epsilon.\epsilon fg.f\eta g \ar[r]^-{\epsilon.Sg} & \epsilon}\end{equation}
in $\AA(B,B)(fg,1)$ and so a unique isomorphism $T:g\epsilon.\eta g\cong 1$
sent by \eqref{e-.f} to \eqref{eq1a}; that is, a unique $T$ satisfying
the first tetrahedron equation.

(iii) Consider the diagram
$$\xymatrix{
&&& gf \ar[dr]^{\eta gf} \ddrruppertwocell<10>^1{Tf} \\
& gf \ar[dr]_{\eta gf} \ar@/^1pc/[urr]^{1} &&& gfgf \ar[dr]^{g\epsilon f} \\ 
1 \ar[ur]^{\eta} \ar[dr]_{\eta} \rrtwocell<\omit>{\eta\eta} && 
gfgf \ar[r]^{gf\eta gf} \urrtwocell<\omit>{~~~gfTf} \ar@/^1pc/[urr]^{1} 
\drrtwocell<\omit>{~~~gSgf} \ar@/_1pc/[drr]_{1} & 
gfgfgf \ar[ur]_{gfg\epsilon f} \rrtwocell<\omit>{~~~g\epsilon\epsilon f} 
\ar[dr]^{g\epsilon fgf} && gf \\
& gf \ar[ur]_{gf\eta} \ar@/_1pc/[drr]_{1}  &&& gfgf \ar[ur]_{g\epsilon f} \\
&&& gf \ar[ur]_{gf\eta} \uurrlowertwocell<-10>_1{gS} }$$
We are to show that the second tetrahedron equation holds. Since $g$ and $f$ are biequivalences, it
will suffice to show that the composite of the central 3 cells above 
is an identity; but since the first tetrahedron equation holds,
this is equivalent to the whole displayed diagram being an
identity. By naturality
of $\eta gf$ and $gf\eta$, this is equal to 
$$
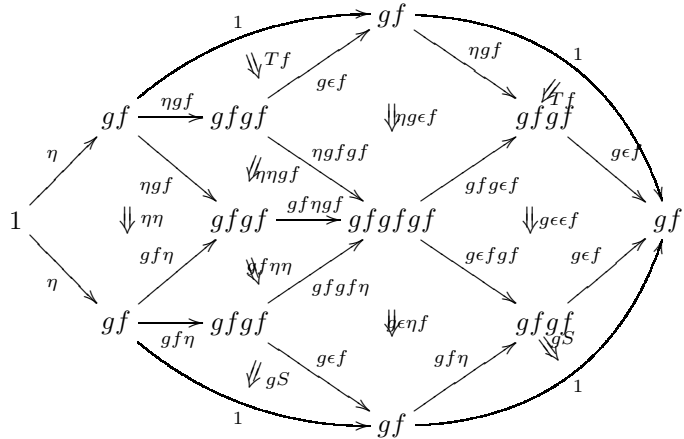
\xymatrix{
&&& gf \ar[dr]^{\eta gf} \ddrruppertwocell<10>^1{Tf} \\
& gf \ar[dr]_{\eta gf} \ar@/^1pc/[urr]^{1} \urrtwocell<\omit>{~~Tf} \ar[r]^{\eta gf} \drrtwocell<\omit>{~~\eta\eta gf} & 
gfgf \ar[ur]_{g\epsilon f} \ar[dr]^{\eta gfgf} \rrtwocell<\omit>{~~\eta g{\epsilon f}} && gfgf \ar[dr]^{g\epsilon f} \\ 
1 \ar[ur]^{\eta} \ar[dr]_{\eta} \rrtwocell<\omit>{~~\eta\eta} && 
gfgf \ar[r]^{gf\eta gf} & 
gfgfgf \ar[ur]_{gfg\epsilon f} \rrtwocell<\omit>{~~~g\epsilon\epsilon f} 
\ar[dr]^{g\epsilon fgf} && gf \\
& gf \ar[ur]^{gf\eta} \ar@/_1pc/[drr]_{1}  \drrtwocell<\omit>{~~gS} \ar[r]_{gf\eta} \urrtwocell<\omit>{gf\eta\eta} & 
gfgf \ar[ur]_{gfgf\eta} \ar[dr]^{g\epsilon f} \rrtwocell<\omit>{g\epsilon{\eta f}}
&& gfgf \ar[ur]^{g\epsilon f} \\
&&& gf \ar[ur]^{gf\eta} \uurrlowertwocell<-10>_1{gS} }$$
which by naturality of the pseudonaturality isomorphisms is equal to 
$$
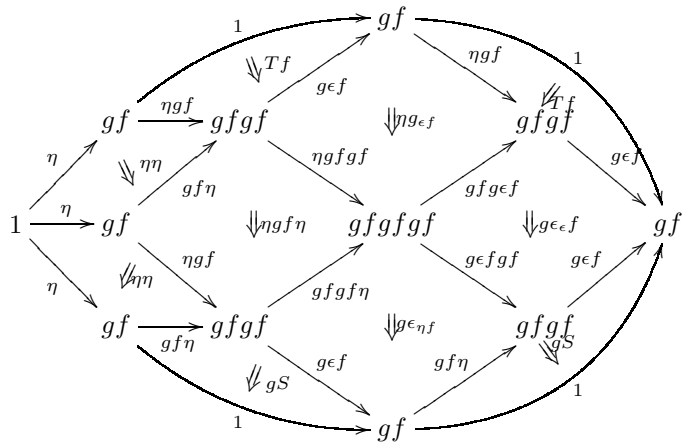
\xymatrix{
&&& gf \ar[dr]^{\eta gf} \ddrruppertwocell<10>^1{Tf} \\
& gf \ar@/^1pc/[urr]^{1} \urrtwocell<\omit>{~~Tf} \ar[r]^{\eta gf} & 
gfgf \ar[ur]_{g\epsilon f} \ar[dr]^{\eta gfgf} \rrtwocell<\omit>{~~\eta g_{\epsilon f}} && gfgf \ar[dr]^{g\epsilon f} \\ 
1 \ar[ur]^{\eta} \ar[dr]_{\eta} \ar[r]^{\eta} \urrtwocell<\omit>{~~\eta\eta} 
\drrtwocell<\omit>{\eta\eta} & gf \ar[ur]_{gf\eta} \ar[dr]^{\eta gf} 
\rrtwocell<\omit>{~~~\eta gf\eta} &&
%gfgf \ar[r]^{gf\eta gf} & 
gfgfgf \ar[ur]_{gfg\epsilon f} \rrtwocell<\omit>{~~~g\epsilon_\epsilon f} 
\ar[dr]^{g\epsilon fgf} && gf \\
& gf \ar@/_1pc/[drr]_{1}  \drrtwocell<\omit>{~~gS} \ar[r]_{gf\eta} & 
gfgf \ar[ur]_{gfgf\eta} \ar[dr]^{g\epsilon f} \rrtwocell<\omit>{~~g\epsilon_{\eta f}}
&& gfgf \ar[ur]^{g\epsilon f} \\
&&& gf \ar[ur]^{gf\eta} \uurrlowertwocell<-10>_1{gS} }$$
and now by naturality of $T$ with respect to $g\epsilon f.gf\eta$ this
is 
$$
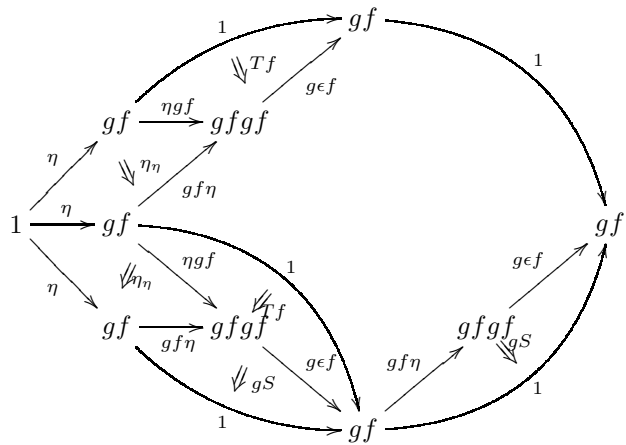
\xymatrix{
&&& gf \ar@/^2pc/[ddrr]^1 \\
& gf \ar@/^1pc/[urr]^{1} \urrtwocell<\omit>{~~Tf} \ar[r]^{\eta gf} & 
gfgf \ar[ur]_{g\epsilon f} \\
1 \ar[ur]^{\eta} \ar[dr]_{\eta} \ar[r]^{\eta} \urrtwocell<\omit>{~~\eta_\eta} 
\drrtwocell<\omit>{\eta_\eta} & gf \ar[ur]_{gf\eta} \ar[dr]^{\eta gf} 
\ddrruppertwocell<10>^1{Tf} &&&& gf \\
& gf \ar@/_1pc/[drr]_{1}  \drrtwocell<\omit>{~~gS} \ar[r]_{gf\eta} & 
gfgf \ar[dr]^{g\epsilon f} && gfgf \ar[ur]^{g\epsilon f} \\
&&& gf \ar[ur]^{gf\eta} \uurrlowertwocell<-10>_1{gS} }$$
which by two applications of the first tetrahedron equation is the identity.
\endproof

\begin{remark}
In fact given $(\eta,S)$ as in (i), and 
another choice $(\eta',S')$, there is a unique invertible $Y:\eta\cong\eta'$
compatible with $S$ and $S'$. Similarly, given $(g',\epsilon')$ in place
of $(g,\epsilon)$, there is a suitable equivalence $g\simeq g'$.
In summary, the ``space of ways of making $f$ into an adjoint biequivalence''
is contractible. We shall not need this fact, and so do not bother to formulate
it precisely.
\end{remark}

We are almost ready to give our condition (ii*). First we record the 
following easy result:

\begin{lemma}\label{lemma:special-lifting}
Let $P:\E\to\B$ be a fibration in \twocat, and $f:D\to E$ an equivalence
in \E. If $\beta:g\cong Pf$ is an invertible 2-cell, and $g$ is part of
an adjoint equivalence $(g:PD\to PE, g^*:PE\to PD, \eta:1\cong g^*g,
\epsilon:gg^*\cong 1)$, then for any invertible lifting
$\overline{\beta}:\overline{g}\cong f$ of $\beta$, we can make
$\overline{g}$ into an adjoint equivalence
$(\overline{g},\overline{g}^*,\eta',\epsilon')$
over $(g,g^*,\eta,\epsilon)$.  
\end{lemma}

\proof
First make the equivalence $f$ into an adjoint equivalence
$(f,f_1,\eta_1,\epsilon_1)$. Then $(Pf,Pf_1,P\eta_1,P\epsilon_1)$ is
an  adjoint equivalence in \BB. The isomorphism $\beta:g\cong Pf$
determines a unique isomorphism $\beta^*:g^*\cong Pf_1$ making
the diagram
$$\xymatrix{
1 \ar[r]^{\eta} \ar[d]_{P\eta_1} & g^* g \ar[d]^{g^*\beta} \\
Pf_1.Pf \ar[r]_{P(\beta^*)^{-1}.Pf} & g^*.Pf }$$
commute. Lift $\beta^*:g^*\cong Pf_1$ to an isomorphism 
$\overline{\beta^*}:\overline{g}^*\cong f_1$ over $\beta^*$. 
Let $\eta':1\to\overline{g}^*\overline{g}$ be the isomorphism
$$\xymatrix{
1 \ar[r]^{\eta_1} & f_1f \ar[r]^{\overline{\beta^*}^{-1}.f} & 
\overline{g}^*.f \ar[r]^{\overline{g}^*.\overline{\beta}^{-1}} & 
\overline{g}^*.\overline{g} }$$
and now observe that $P\eta'=\eta$ by the defining property of 
$\beta^*$.

Now $\overline{g}$ is isomorphic to the equivalence $f$, and so is
itself an equivalence. We have an isomorphism
$\eta':1\cong\overline{g^*}\overline{g}$,
and so there is a unique isomorphism
$\epsilon':\overline{g}\overline{g^*}\cong 1$ for which 
$(\overline{g},\overline{g^*},\eta',\epsilon')$ is an adjoint
equivalence in \EE. But then $(P\overline{g},P\overline{g^*},P\eta',P\epsilon')$
and $(g,g^*,\eta,\epsilon)$ are both adjoint equivalences in \BB,
with $P\overline{g}=g$, $P\overline{g^*}=g^*$, and $P\eta'=\eta$; 
so finally $P\epsilon'=\epsilon$, since the unit of an adjunction determines the counit.
\endproof

Our condition (ii*) is now that $F:\AA\to\BB$ has the 
{\em adjoint-biequivalence-lifting property\/}: given $E\in\EE$, $B\in\BB$,
$f:B\to FE$ and $g:FE\to B$ along with $\eta$, $\epsilon$, $S$, and
$T$ giving an adjoint biequivalence in \BB, there exist 
$D\in\EE$ and $f':D\to E$, $g':E\to D$, $\eta':1\to g'f'$,
$\epsilon':f'g'\to 1$, $S':\epsilon'f'.f'\eta'$, $T':1\to g'\epsilon'.\eta' g'$
forming an adjoint biequivalence in \EE, with $Ff'=f$, $Fg'=g$, 
$F\eta'=\eta$, $F\epsilon'=\epsilon$, $FS'=S$, and $FT'=T$. 

Any \Gray-functor $F:\AA\to\BB$ which is a fibration on the homs and has
the adjoint-biequivalence-lifting property is certainly a fibration:
for given $A\in\AA$ and an isomorphism $B\cong FA$ in $\pi_*\BB$, 
we can represent the isomorphism by a biequivalence $f:B\to FA$ in \BB,
and then extend this to an adjoint biequivalence (involving $g$, $\eta$,
$\epsilon$, and so on), which by assumption can be lifted to an adjoint
biequivalence involving $f':A'\to A$ and other data, lying above the 
original adjoint biequivalence. In particular $f'$ represents an
isomorphism $A'\cong A$ in $\pi_*\AA$ lying over the original 
isomorphism $B\cong FA$ in $\pi_*\BB$. The harder part is:

\begin{proposition}
Every fibration satisfies the adjoint-biequivalence-lifting property.
\end{proposition}

\proof
Let $P:\EE\to\BB$ be a fibration.
Let $E$, $B$, $f$, $g$, $\eta$, $\epsilon$, $S$, and $T$ be an adjoint
biequivalence as above. Then $f$ determines an isomorphism
$B\cong FE$ in $\pi_*\BB$, which can be lifted to an isomorphism 
$D\cong E$ in $\pi_*\EE$, and this isomorphism can be represented
by a biequivalence $f_0:D\to E$ in \EE, with $PD=B$ and $Pf_0$
equivalent in $\BB(B,PE)$ to $f$. Since $P:\EE(D,E)\to\BB(B,PE)$
is a fibration in \twocat, we can lift the equivalence $Pf_0\simeq f$
to an equivalence $f_0\simeq f'$, so that $Pf'=f$; also $f'$ is
equivalent  to the biequivalence $f_0$ so is itself a biequivalence.
Thus we have lifted the biequivalence $f:B\to PE$ to a biequivalence
$f':D\to E$; we must now lift the remaining components of the adjoint 
biequivalence.

By Proposition~\ref{prop:adj-biequiv},
 we can make $f'$ into an adjoint biequivalence in \EE via $(g_1:E\to
D, \eta_1, \epsilon_1, S_1, T_1)$, and now its image $(f,Pg_1,P\eta_1,
P\epsilon_1, PS_1, PT_1)$ under $P$ is an adjoint biequivalence in
\BB. But $(f,g,\eta,\epsilon,S,T)$ is also an adjoint biequivalence, so the composite 
$$\xymatrix{
g \ar[r]^-{P\eta_1.g} & Pg_1.f.g \ar[r]^-{Pg_1.\epsilon} & Pg_1 }$$
is (part of) an adjoint equivalence, and so can be lifted to an 
adjoint equivalence $\alpha:g'\to g_1$ in $\EE(D,E)$, which in turn gives
an adjoint equivalence
$$\xymatrix{
f'g' \ar[r]^{f'\alpha} & f'g_1 \ar[r]^{\epsilon_1} & 1 }$$
lying over the top leg of 
$$\xymatrix{
Pf'.Pg' \ar@{=}[d] \ar[drr]^{Pf'.P\alpha} \\
fg \ar[r]_-{f.P\eta'.g} \drlowertwocell_{1}{PT'.g} & f.Pg_1.f.g \ar[r]_-{f.Pg_1.\epsilon} \ar[d]_{P\epsilon.fg}  \drlowertwocell<\omit>{} &
f.Pg_1 \ar[d]^{P\epsilon_1} \\
& fg \ar[r]_{\epsilon} & 1
}$$
in which the unnamed 2-cell is a pseudonaturality isomorphism. The pasting composite
here is an invertible 2-cell in $\BB(B,B)$ between adjoint equivalences, so can be 
lifted to an invertible 2-cell 
$$\xymatrix{
f'g' \ar[r]^{f'\alpha} \rrlowertwocell_{\epsilon'}{X} & 
f'g_1 \ar[r]^{\epsilon_1} & 1 }$$
in $\EE(D,D)$, and now the adjoint equivalence structure of 
$\epsilon_1.f'\alpha$ transports across the isomorphism $X$ to give
an adjoint equivalence structure on $\epsilon'$, lying over that
on $\epsilon:fg\to 1$, as in Lemma~\ref{lemma:special-lifting}.

At this point we have $f':D\to E$ and  $g':E\to D$ over $f$ and $g$,
and an adjoint equivalence $\epsilon':f'g'\to 1$ over $\epsilon:fg\to 1$.
Now $f'$ is a biequivalence, so by Proposition~\ref{prop:adj-biequiv}
there exist an equivalence $\beta:1\to g'f'$ and an invertible 2-cell
$$\xymatrix{
f' \ar[r]^{f'\beta} \drlowertwocell_1{S_2} & f'g'f' \ar[d]^{\epsilon'f'} \\
& f' }$$
which $P$ sends to the left hand side of the diagram 
$$\xymatrix{
f \ar[r]^{f.P\beta} \drlowertwocell_1{~~PS_2} & fgf \ar[d]^{\epsilon f} &
f \ar[r]^{f\eta} \drlowertwocell_1{S} & fgf \ar[d]^{\epsilon f} \\
& f & & f }$$
but since we also have the diagram on the right, and $\epsilon$ is an
equivalence and $f$ a biequivalence, there exists a unique invertible
$$\xymatrix{
1 \rtwocell^{\eta}_{P\beta}{R} & gf }$$
which when pasted to $PS_2$ gives $S$. We can now lift this by
Lemma~\ref{lemma:special-lifting} to 
$$\xymatrix{
1 \rtwocell^{\eta'}_{\beta}{R'} & g'f' }$$
with $\eta'$ an adjoint equivalence over $\eta$, and define 
$S'$ to be the composite
$$\xymatrix @C4pc {
f' \rtwocell^{f'\eta'}_{f'\beta}{~~f'R'} \drlowertwocell_1{S_2} &
f'g'f' \ar[d]^{\epsilon'f'} \\
& f' }$$
and observe that this lies over $S$. 

By Proposition~\ref{prop:adj-biequiv}  once again there is a unique invertible
$T':g'\epsilon'.\eta'g'\cong 1$ for which 
$(f',g',\eta',\epsilon',S',T')$ is an adjoint biequivalence; thus 
$(Pf'=f,Pg'=g, P\eta'=\eta, P\epsilon'=\epsilon, PS'=S, PT')$ is
an adjoint biequivalence, and by the uniqueness part of the 
last clause of Proposition~\ref{prop:adj-biequiv}, it follows that $PT'=T$.
\endproof

This now proves:

\begin{proposition}
A \Gray-functor $F:\AA\to\BB$ is a fibration if and only if it satisfies
the adjoint-biequivalence-lifting property and each
$F:\AA(A,B)\to\BB(FA,FB)$ is a fibration in \twocat.
\end{proposition}

We have already observed that the condition on the homs is a right 
lifting property. For the adjoint biequivalence lifting property, note that 
the structure of adjoint biequivalence can be described in terms of
objects, 1-cells, 2-cells, and 3-cells of \AA satisfying certain equations;
thus since \graycat is locally finitely presentable, there is a 
``universal adjoint biequivalence'': a \Gray-category \EE, with the property
that adjoint biequivalences in \AA are in natural bijection with 
\Gray-functors from \EE to \AA. There are two objects of \EE; the 
adjoint-biequivalence-lifting property is exactly the right lifting 
property with respect to one (either) of the \Gray-functors $1\to\EE$.
By the small object argument we conclude:

\begin{proposition}
There is a cofibrantly generated weak factorization system on \graycat 
whose right part is the fibrations. (The left part is what we call 
the trivial cofibrations.)
\end{proposition}

\begin{remark}
Since all objects in \twocat are fibrant, and since there are no 
non-trivial  biequivalences in the terminal \Gray-category, it is 
clear that all objects in \graycat are also fibrant.  
\end{remark}

\section{Proof of the model structure}
\label{sect:proof}

It remains to show that the trivial cofibrations are precisely those
cofibrations which are also weak equivalences. Since every trivial 
fibration is a fibration, certainly every trivial cofibration is 
a cofibration. We shall show that every trivial cofibration is a weak
equivalence, and then use a standard argument to show that also 
every weak equivalence which is a cofibration is a trivial cofibration.

The key step in this proof is the existence of path objects:

\begin{proposition}
Path objects exist in \graycat: for every \Gray-category \BB there exists
a \Gray-category \PB and a factorization
$$\xymatrix{
\BB \ar[r]^{D} & {}\PB \ar[r]^-{\binom{P}{P'}} & \BB\times\BB }$$
of the diagonal into a weak equivalence followed by a fibration.
\end{proposition}

\proof
Define a \Gray-category \PB as follows:
\begin{itemize}
\item an object is a biequivalence $a:A\to A'$ in \BB
\item a 1-cell from $a:A\to A'$ to $b:B\to B'$ consists of 1-cells
$f:A\to B$ and $f':A'\to B'$ in \BB, and an equivalence 
$\phi:bf\simeq f'a$
\item a 2-cell from $(f,f',\phi)$ to $(g,g',\psi)$ consists of 2-cells
$\xi:f\to g$ and $\xi':f'\to g'$ equipped with an invertible 3-cell $\Xi$
between
$$\xymatrix {
A \ar[r]^{a} \dtwocell_f^g{^\xi} & A' \duppertwocell^{g'}{^\psi} && 
A \ar[r]^{a} \dlowertwocell_{f}{^\phi} & 
A' \dtwocell_{f'}^{g'}{^\xi'} \\
B \ar[r]_{b} & B' && B \ar[r]_{b} & B' }$$
\item a 3-cell from $(\xi,\xi',\Xi)$ to $(\zeta,\zeta',Z)$ is a
pair of 3-cells $M:\xi\to\zeta$ and $M':\xi'\to\zeta'$ satisfying 
the evident compatibility condition.
\end{itemize}
This is made into a \Gray-category in the obvious way. One now checks that
\begin{itemize}
\item a 1-cell $(f,f',\phi)$ in \PB is a biequivalence if and only if
$f$ and $f'$ are biequivalences in \B
\item a 2-cell $(\xi,\xi',\Xi)$ in \PB is an equivalence if and only if
$\xi$ and $\xi'$ are equivalences in \B
\item a 3-cell $(M,M')$ in \PB is invertible if and only if $M$ and $M'$
are so in \B.
\end{itemize}
There are 
evident projections $P,P':\PB\to\BB$ which are \Gray-functors, and the 
diagonal $\Delta:\BB\to\BB\times\BB$ factorizes as 
$$\xymatrix{
\BB \ar[r]^{D} & {}\PB \ar[r]^{\binom{P}{P'}} & \BB\times\BB }$$
It follows easily from the characterization of biequivalences, equivalences,
and isomorphisms in \PB given above, that $\binom{P}{P'}$ is a fibration.
It is also straightforward to check that $D$ is a weak equivalence.
\endproof

We now complete the proof of the theorem using a standard argument.
We are to show that
the trivial cofibrations are precisely the cofibrations which are 
weak equivalences.

Suppose first that $F$ is a trivial cofibration; we have already observed
that $F$ is a cofibration, so we must show that it is also a weak equivalence.
Since \AA (like any object) is fibrant, there exists a map $G:\BB\to\AA$
with $GF=1$. 

Now $F$ is a trivial cofibration and $\binom{P}{P'}$ a fibration, so there
is a lifting $H$ as in 
$$\xymatrix{
\AA \ar[d]_F \ar[r]^F & \BB \ar[r]^D & \PB \ar[d]^{\binom{P}{P'}} \\
\BB \ar@{.>}[urr]^H \ar[rr]_{\binom{1}{FG}} && \BB\times\BB. }$$
Since $D$ is a weak equivalence and $PD=P'D=1$, both $P$ and $P'$
are weak equivalences. Since $PH=1$ and $P$ is a weak equivalence,
$H$ is a weak equivalence. Since $P'H=FG$ and $P'$ and $H$ are weak
equivalences, $FG$ is a weak equivalence. 

Finally $F$ is a retract of $FG$, since $GF=1$, and so $F$ is indeed
a weak equivalence. This proves that every trivial cofibration is a 
weak equivalence and a cofibration; we now turn to the converse.

If $F$ is both a weak equivalence and a cofibration, factorize it
as $F=QJ$, where $J$ is a trivial cofibration and $Q$ a fibration. 
Now $F$ and $J$ are weak equivalences, so $Q$ is a weak equivalence,
and hence a trivial fibration. But now since $F$ is a cofibration, and
factorizes through $Q$ by $J$, it is a retract of $J$ and hence is itself
a trivial cofibration.

This completes the proof of Theorem~\ref{thm:main}.
\endproof

\section{\Gray-groupoids}
\label{sect:groupoids}

In this section we restrict from general \Gray-categories to \Gray-groupoids:
these are \Gray-categories in which all l-cells, 2-cells, and 3-cells are
(strictly) invertible. The \Gray-groupoids form a full subcategory 
\graygpd of \graycat. Since the structure of of \Gray-category and that of
\Gray-groupoid are both essentially algebraic, and the inclusion is given by forgetting 
some of this algebraic structure, it has a left adjoint, and so the full
subcategory is reflective. We shall show that the model structure on 
\graycat restricts to \graygpd.

\Gray-groupoids were called {\em algebraic homotopy 3-types} in 
\cite{GPS}, and we shall see that they do indeed provide a model
for homotopy 3-types --- this is due originally to Joyal and Tierney,
in unpublished but widely advertised work; here, we shall give a
 functorial and model-theoretic proof. The relationship between 
\Gray-groupoids and homotopy 3-types has also been studied
in \cite{Berger-3type} and \cite{Leroy}; see the introduction for comments 
about the relationship with \cite{Berger-3type}.

\begin{theorem}
There is a combinatorial model structure on \graygpd for which a morphism 
is a fibration or weak equivalence if and only if it is one in \graycat.
\end{theorem}

\proof
The cofibrations are then morphisms with the left lifting property with
respect to the trivial fibrations, and the trivial cofibrations are
the morphisms with the left lifting property with respect to the
fibrations. The existence of the weak factorization systems is
immediate. Once again every trivial cofibration
is a cofibration; we have to show that a cofibration is trivial if and only
if it is a weak equivalence. This will follow exactly as before provided
that we can construct path objects in \graygpd; and this will certainly 
be the case if our path object \PB (in \graycat) is a \Gray-groupoid whenever \BB is one.

Suppose then that \BB is a \Gray-groupoid, and consider the \Gray-category
\PB. A 3-cell (between specified 2-cells) consists of a pair $(M,M')$ 
of 3-cells in \BB satisfying a compatibility condition. Since \BB is 
a \Gray-groupoid, both $M$ and $M'$ are invertible, and so the 3-cell
$(M,M')$ in \PB is invertible. A 2-cell, between specified 1-cells
$(f,\phi):a\to b$ and $(f',\phi':a\to b)$, consists
of 2-cells $\xi:f\to g$ and $\xi':f'\to g'$ in \BB equipped with an invertible
3-cell $\Xi$. Since \BB is a \Gray-groupoid, $\xi$ and $\xi'$ have inverses
$\xi^{-1}$ and $\xi'^{-1}$, and these become an inverse to 
$(\xi,\xi',\Xi)$ when equipped with the 3-cell
$\xi'^{-1}a.\Xi^{-1}.b\xi^{-1}:\phi.b\xi^{-1}\cong \xi'^{-1}a.\psi$. Thus all
2-cells in \PB are invertible. It remains to show that every 1-cell in
\PB is invertible. A 1-cell from $a:A\to A'$ to $b:B\to B'$ consists of 
1-cells $f:A\to B$ and $f':A'\to B'$ in \BB, equipped with an equivalence
$\phi:bf\simeq f'a$. But since \BB is a \Gray-groupoid, $f$ and $f'$
have inverses $f^{-1}$ and $f'^{-1}$, and $\phi$ is not just an equivalence
but an isomorphism; and now $f^{-1}$ and $f'^{-1}$ become an inverse 
for $(f,f',\phi)$ when equipped with the 2-cell
$f'^{-1}.\phi^{-1}.f^{-1}:af^{-1}=f'^{-1}f'af^{-1}\cong f'^{-1}bff^{-1}=f'^{-1}b$.
\endproof

In this restricted version, the description of fibrations becomes
simpler. A \Gray-functor $P:\EE\to\BB$ between \Gray-groupoids is a fibration
when the following conditions hold:
\begin{itemize}
\item for every $E\in\EE$ and every 1-cell $f:B\to PE$ in \BB there is a lifting
$f':D\to E$ in \EE
\item for every $f:D\to E$ in \EE and every 2-cell $\beta:b\to Pf$ in \BB
there is a lifting $\gamma:d\to f$ in \EE
\item for every 2-cell $\alpha:f\to g:D\to E$ in \EE and every 3-cell
$M:\beta\to P\alpha$ in \BB there is a lifting $M':\gamma\to\alpha$
in \EE
\end{itemize}
A \Gray-functor $F:\EE\to\BB$ between \Gray-groupoids is a weak equivalence
when the following conditions hold:
\begin{itemize}
\item for every $B\in\BB$ there is an $E\in\EE$ and a morphism $f:B\to FE$
\item for all $D,E\in\EE$ and every $b:FD\to FE$ in \BB there is a 
morphism $f:D\to E$ in \EE and a 2-cell $\beta:b\to Pf$ in \BB
\item for all $\alpha,\beta:f\to g:D\to E$ in \EE and every 3-cell 
$Y:P\alpha\to P\beta$ in \BB there is a 3-cell $X:\alpha\to\beta$ with
$PX=Y$.
\end{itemize}

In \cite{Berger-3type}, an adjunction between Gray-groupoids and simplicial
sets was described, with the right adjoint $N:\graygpd\to\SSet$ giving
the nerve of a \Gray-groupoid; we shall write $\Pi_3$ for the left adjoint. 

\begin{proposition}
The nerve functor $N:\graygpd\to\SSet$ preserves fibrations and trivial
fibrations, and so is the right adjoint part of a Quillen adjunction.
\end{proposition}

\proof
Suppose first that $F:\AA\to\BB$ is a trivial fibration in \graygpd. We
must show that $NF:N\AA\to N\BB$ is a trivial fibration of simplicial
sets; in other words, that it has the right lifting property with respect
to the inclusion $\partial\Delta[n]\to\Delta[n]$ for all $n$. For $n=0$
this is the fact that $F$ is surjective on objects and for $n=1$ it is
the fact that $F$ is full on 1-cells. For $n=2$ it says that for any 
1-cells $f:A\to B$, $g:B\to C$ and $h:A\to C$ in \AA, and any 2-cell
$\beta:Fg.Ff\to Fh$ in \BB, there is a 2-cell $\alpha:gf\to h$ in \AA
with $F\alpha=\beta$; this holds since $F$ is full on 2-cells ($F$ is 
locally full). For $n=3$ it becomes a little more complicated. Consider
2-cells in \AA as in the diagram
% $$
% \xy
% (0,-8)*{A_0}="A0";
% (12,8)*{A_1}="A1";
% (28,8)*{A_2}="A2";
% (40,-8)*{A_3}="A3";
% {\ar^{f_{01}} "A0";"A1"};
% {\ar^{f_{12}} "A1";"A2"};
% {\ar^{f_{23}} "A2";"A3"};
% {\ar_{f_{03}} "A0";"A3"};
% {\ar_{f_{02}} "A0";"A2"};
% {\ar@{=>}_{\beta} (28,2);(28,-2)};
% {\ar@{=>}_{\alpha} (15,6);(15,2)};
% \endxy
% \qquad %= \qquad
% \xy
% (0,-8)*{A_0}="A0";
% (12,8)*{A_1}="A1";
% (28,8)*{A_2}="A2";
% (40,-8)*{A_3}="A3";
% {\ar^{f_{01}} "A0";"A1"};
% {\ar^{f_{12}} "A1";"A2"};
% {\ar^{f_{23}} "A2";"A3"};
% {\ar_{f_{03}} "A0";"A3"};
% {\ar_{f_{13}} "A1";"A3"};
% {\ar@{=>}^{\gamma} (12,2);(12,-2)};
% {\ar@{=>}^{\delta} (25,6);(25,2)};
% \endxy
% $$

$$
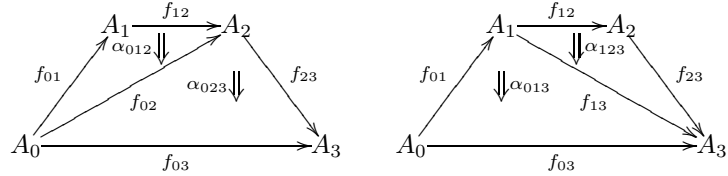
\xy
(0,-8)*{A_0}="A0";
(12,8)*{A_1}="A1";
(28,8)*{A_2}="A2";
(40,-8)*{A_3}="A3";
{\ar^{f_{01}} "A0";"A1"};
{\ar^{f_{12}} "A1";"A2"};
{\ar^{f_{23}} "A2";"A3"};
{\ar_{f_{03}} "A0";"A3"};
{\ar_{f_{02}} "A0";"A2"};
{\ar@{=>}_{\alpha_{023}} (28,2);(28,-2)};
{\ar@{=>}_{\alpha_{012}} (18,7);(18,3)};
\endxy
\qquad %= \qquad
\xy
(0,-8)*{A_0}="A0";
(12,8)*{A_1}="A1";
(28,8)*{A_2}="A2";
(40,-8)*{A_3}="A3";
{\ar^{f_{01}} "A0";"A1"};
{\ar^{f_{12}} "A1";"A2"};
{\ar^{f_{23}} "A2";"A3"};
{\ar_{f_{03}} "A0";"A3"};
{\ar_{f_{13}} "A1";"A3"};
{\ar@{=>}^{\alpha_{013}} (12,2);(12,-2)};
{\ar@{=>}^{\alpha_{123}} (22,7);(22,3)};
\endxy
$$
% $$\xymatrix{
% & A_1 \ar[r]^{f_{12}} & A_2 \ar[dr]^{f_{23}} &&
% & A_1 \ar[r]^{f_{12}} \ar[drr]_{f_{13}} \drrtwocell\omit{<-2>\alpha_{123}} & A_2 \ar[dr]^{f_{23}} & \\
% A_0 \ar[ur]^{f_{01}} \ar[urr]_{f_{02}} \ar[rrr]_{f_{03}} \urrtwocell\omit{<-2>\alpha_{012}} & {}\rtwocell\omit{<-2>\alpha_{023}} &{} & A_3 &
% A_0 \ar[ur]^{f_{01}} \ar[rrr]_{f_{03}} &
% {}\rtwocell\omit{<-2>\alpha_{013}} &{} & A_3 }$$
where the left composite may be written as $\alpha_{023}.f_{23}\alpha_{012}$ 
and the right composite as $\alpha_{013}.\alpha_{123} f_{01}$. Apply $F$
to each side, and suppose that we have a 3-cell 
$Y:F(\alpha_{023}.f_{23}\alpha_{012})\to F(\alpha_{013}.\alpha_{123} f_{01})$
between the resulting 2-cells in \BB. The case $n=3$ amounts to the fact
that there is a 3-cell $X:\alpha_{023}.f_{23}\alpha_{012}\to
\alpha_{013}.\alpha_{123} f_{01}$ 
in \AA with $FX=Y$; this is true since
$F:\AA\to\BB$ is full on 3-cells (locally locally full). For $n>3$, an
$n$-simplex in the nerve of a \Gray-groupoid does not involve further data,
it just amounts to the assertion that two pasting composites (with the same 
2-cells as domain and codomain) of 3-cells
are equal. Since $F$ is faithful on 3-cells (locally locally faithful),
such assertions can be lifted from \BB to \AA, and so all remaining conditions
hold. This proves that $NF$ is a trivial fibration.

Suppose now that $F:\AA\to\BB$ is a fibration in \graygpd. We must show 
that $NF:N\AA\to N\BB$ is a fibration of simplicial sets (a Kan fibration);
in other words, that it has the left lifting property with respect to the
horns $\Lambda^r[n]\to\Delta[n]$.

For $n=1$, this says that given $A\in\AA$ and a 1-cell $f:B\to FA$ in \BB
there is a 1-cell $f':A'\to\AA$ with $FA'=B$ and $Ff'=f$, as well as 
a corresponding statement involving 1-cells $FA\to B$. This first is part 
of our characterization of fibrations; the second an easy consequence
of it, given that every 1-cell is invertible.

For $n=2$ and $r=1$, this says that given $f:A\to B$ and $g:B\to C$ in \AA
and a 2-cell $\beta:Fg.Ff\to k$ in \BB there is a 1-cell $h:A\to C$
with $Fh=k$ and a 2-cell $\alpha:gf\to h$ with $F\alpha=\beta$. This is
once again part of our characterization of fibrations in \graygpd.

For $n=2$ and $r=0$ we start with 1-cells $g:A\to B$ and $u:A\to C$ in \AA,
and a 1-cell $v:FB\to FC$ and 2-cell $\gamma:v.Ff\to Fu$ in \BB, as in
the solid parts of the diagram below. 
$$
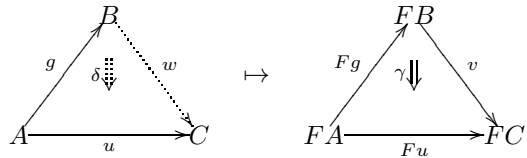
\xy
(0,-8)*{A}="A";
(12,8)*{B}="B";
(24,-8)*{C}="C";
{\ar^{g} "A";"B"};
{\ar@{.>}^{w} "B";"C"};
{\ar_{u} "A";"C"};
{\ar@{:>}_{\delta} (12,2);(12,-2)};
\endxy
\quad\mapsto\quad
\xy
(0,-8)*{FA}="A";
(12,8)*{FB}="B";
(24,-8)*{FC}="C";
{\ar^{Fg} "A";"B"};
{\ar^{v} "B";"C"};
{\ar_{Fu} "A";"C"};
{\ar@{=>}_{\gamma} (12,2);(12,-2)};
\endxy
$$
We must
construct a 1-cell $w:B\to C$ in \AA with $Fw=v$ and a 2-cell
$\delta:wg\to u$ with $F\delta=\gamma$.
Let 
$f=u^{-1}$ and $k=v^{-1}$. Then $v^{_1}.v.Fg.Ff=Fg.Ff$ and 
$v^{-1}.Fu.Ff=k.Fu.Fu^{-1}=k$, and so the 2-cell $v^{-1}.\gamma.Ff$ goes
from $Fg.Ff$ to $k$. By the previous paragraph there is a 2-cell
$\alpha:gf\to h$ in \AA with $Fh=k$ and $F\alpha=v^{-1}.\gamma.Ff$. 
Let $w=h^{-1}$ and $\delta=w\alpha u:wgfu\to whu$. Then $wgfu=wgu^{-1}u=wg$ and 
$whu=h^{-1}hu=u$, and so $\delta$ is indeed a 2-cell from $wg$ to $u$. Also $Fw=Fh^{-1}=(Fh)^{-1}=k^{-1}=v$, while $F\delta=Fw.F\alpha.Fu=Fw.v^{-1}.\gamma.Ff.Fu=Fh^{-1}.k.\gamma.Fu^{-1}.Fu=\gamma$,
and so the condition does hold. The case $n=r=2$ is similar.

For $n=3$, we start with objects $A_0$, $A_1$, $A_2$, and $A_3$
of \AA, with 1-cells $f_{ij}:A_i\to A_j$ for $0\le i<j\le 3$. For
$0\le i<j<k\le 3$ with $r\in\{i,j,k\}$ we have a 2-cell 
$\alpha_{ijk}:f_{jk}f_{ij}\to f_{ik}$ in \AA. For the choice of 
$0\le i<j<k\le 3$ with $r\notin\{i,j,k\}$ there is a 2-cell 
$\beta_{ijk}:F(f_{jk}f_{ij})\to Ff_{ik}$ in \BB, and finally a 3-cell $Y$. 
between the two pasting composites.

We treat the case $r=1$ in detail. The ``missing 2-cell'' is 
$\alpha_{023}:f_{23}f_{02}\to f_{03}$. Since $\alpha_{012}$ is an 
equivalence, we can find a 2-cell $\alpha'_{023}:f_{23}f_{02}\to f_{03}$
and a 3-cell 
$X':\alpha'_{023}.f_{23}\alpha_{012}\to \alpha_{013}.\alpha_{123} f_{01}$.
Applying $F$ gives a 3-cell
$FX':F\alpha'_{023}.F(f_{23}\alpha_{012})\to F\alpha_{013}.F(\alpha_{123} f_{01})$,
but we also have a 3-cell
$Y:\beta_{023}.F(f_{23}\alpha_{012})\to F\alpha_{013}.F(\alpha_{123} f_{01})$, and
now there is a unique 3-cell $Z:\beta_{023}\to F\alpha'_{023}$ which pastes
onto $FX'$ to give $Y$. We can lift this to a 3-cell 
$W:\alpha_{023}\to\alpha'_{023}$ with $F\alpha_{023}=\beta_{023}$ and
$FW=Z$. Finally pasting $W$ onto $X'$ gives a 3-cell
$X:\alpha_{023}.f_{23}\alpha_{012}\to \alpha_{013}.\alpha_{123} f_{01}$ 
over $Y$.

The case $r=2$ is similar to that of $r=1$. For $r=0$, where $\alpha_{123}$ is
the missing face, we can take inverses and work instead with 
$$
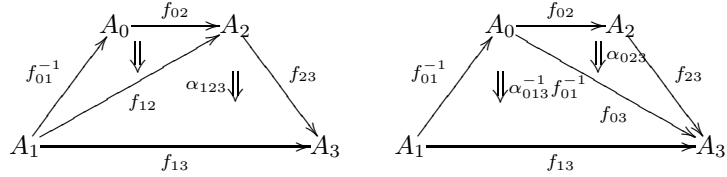
\xy
(0,-8)*{A_1}="A0";
(12,8)*{A_0}="A1";
(28,8)*{A_2}="A2";
(40,-8)*{A_3}="A3";
{\ar^{f^{-1}_{01}} "A0";"A1"};
{\ar^{f_{02}} "A1";"A2"};
{\ar^{f_{23}} "A2";"A3"};
{\ar_{f_{13}} "A0";"A3"};
{\ar_{f_{12}} "A0";"A2"};
{\ar@{=>}_{\alpha_{123}} (28,2);(28,-2)};
{\ar@{=>} (15,6);(15,2)};
\endxy
\qquad %= \qquad
\xy
(0,-8)*{A_1}="A0";
(12,8)*{A_0}="A1";
(28,8)*{A_2}="A2";
(40,-8)*{A_3}="A3";
{\ar^{f^{-1}_{01}} "A0";"A1"};
{\ar^{f_{02}} "A1";"A2"};
{\ar^{f_{23}} "A2";"A3"};
{\ar_{f_{13}} "A0";"A3"};
{\ar_(0.6){f_{03}} "A1";"A3"};
{\ar@{=>}^{\alpha^{-1}_{013}f^{-1}_{01}} (12,2);(12,-2)};
{\ar@{=>}^{\alpha_{023}} (25,6);(25,2)};
\endxy
$$
where the unnamed face on the left is $\alpha^{-1}_{012}f^{-1}_{01}$;
% $$\xymatrix{
% & A_0 \ar[r]^{f_{02}} & A_2 \ar[dr]^{f_{23}} &&
% & A_0 \ar[r]^{f_{02}} \ar[drr]_{f_{03}} \drrtwocell\omit{<-2>\alpha_{023}} & A_2 \ar[dr]^{f_{23}} & \\
% A_1 \ar[ur]^{f^{-1}_{01}} \ar[urr]_{f_{12}} \ar[rrr]_{f_{13}} \urrtwocell\omit{<-2>\alpha^{-1}_{012}f^{-1}_{01}} &  &{} & A_3 &
% A_1 \ar[ur]^{f^{-1}_{01}} \ar[rrr]_{f_{13}} &
% {}\rtwocell\omit{<-2>\alpha_{013}} &{} & A_3 }$$
which reduces it to the case $r=1$. Similarly the case $r=3$ can be 
reduced to $r=2$ by taking inverses. 

Since we are dealing with nerves of \Gray-groupoids, all horns have unique fillers for $n\ge 4$, and so the condition is automatic.
\endproof

By Ken Brown's lemma \cite[Lemma~1.1.12]{Hovey-book}, $N$ takes all weak equivalences between
fibrant objects to weak equivalences; but all objects of \graygpd are
fibrant, and so in fact $N$ also preserves weak equivalences. By a similar
argument, using the fact that all objects of \SSet are cofibrant, $\Pi_3$ also
preserves weak equivalences. We record this as

\begin{proposition}
Both the nerve functor $N:\graygpd\to\SSet$ and its left adjoint
$\Pi_3:\SSet\to\graygpd$ preserve weak equivalences.  \endproof
\end{proposition}

Quillen adjunctions induce derived adjunctions between the homotopy 
categories. Since $\Pi_3$ preserves weak equivalences and all
objects of \SSet are cofibrant,
the simplicial sets $X$ for which the unit of the derived adjunction is
invertible are precisely those for which the unit $X\to N\Pi_3X$ of the 
Quillen adjunction is a weak equivalence. Berger showed in 
\cite[Proposition~3.2]{Berger-3type} that for a simplicial set $X$,
the unit $X\to N\Pi_3X$ is a weak equivalence if and only if $X$ is a 
3-type.

Similarly, since $N$ preserves weak equivalences and all objects of \graygpd are fibrant, the \Gray-groupoids \AA for which the counit of the derived
adjunction is invertible are precisely those for which the counit
$\Pi_3N\AA\to\AA$ is a weak equivalence.
Once we have shown that this holds for all \Gray-groupoids, then we shall have constructed
an equivalence between the homotopy categories of \Gray-groupoids and 
the homotopy category of simplicial 3-types. 

\begin{theorem}
The Quillen adjunction between \graygpd and \SSet induces an equivalence
between the homotopy categories of algebraic homotopy 3-types and 
simplicial homotopy 3-types.
\end{theorem}

\proof 
We need to show that the counit map $E:\Pi_3N\AA\to\AA$ is a triequivalence
of \Gray-categories. 
An object of $\Pi_3N\AA$ is a 0-simplex of $N\AA$; that is, an object of \AA. 
Thus $\Pi_3N\AA$ and \AA have the same objects, and the counit map 
$E:\Pi_3N\AA\to\AA$ is the identity on objects.

The 1-cells of $\Pi_3N\AA$ are freely generated by the non-degenerate 1-simplices
of $N\AA$ together with formal inverses; that is, by the non-identity 1-cells of \AA together with formal inverses. We write a typical such 1-cell
as $\bb f$ or $[f_n]\ldots[f_1]$, where the $f_i$ are 1-cells in \AA or
their formal inverses. 
Any non-identity 1-cell in \AA can be realized as the composite of such a word of length 1
(or 0 in the case of an 
identity), and so the counit $E:\Pi_3N\AA\to\AA$ is surjective (and full) on 
1-cells.

% It remains to show that $E:\Pi_3N\AA\to\AA$ is full on 2-cells, and fully
% faithful on 3-cells. To do this, it is useful to observe that for
% the \Gray-functor $E:\Pi_3N\AA\to\AA$, if 
% $\alpha_1$ and $\alpha_2$ are 2-cells in $\Pi_3N\AA$, with the same domain
% and codomain 1-cells, and $\alpha_1\cong\alpha'_1$ and $\alpha_2\cong\alpha'_2$
% in $\Pi_3N\AA$, then $F$ is full on 2-cells $\alpha_1\to\alpha_2$ if and only 
% if it is full on 2-cells $\alpha'_1\to\alpha'_2$; with a corresponding
% result for faithfulness. 

The 2-cells of $\Pi_3N\AA$ are obtained by
freely adjoining an invertible 2-cell $[f_2][f_1]\to [f]$ for each
2-cell $\phi:f_2f_1\to f$ which is ``non-degenerate'' in the sense
that no more than one of $f_1$, $f_2$, and $\phi$ are identities.

Now a general 1-cell in $\Pi_3N\AA$ is a word $[f_n]\ldots[f_2][f_1]$ in 
the 1-cells of \AA and their formal inverses. But in fact any such
1-cell is isomorphic in $\Pi_3N\AA$ to an identity or to a word consisting
of a single 1-cell in \AA. To see this, first observe that if 
$f:A\to B$ is any 1-cell in \AA, then it has an inverse, say $g$,
and now the equations $gf=1$ and $fg=1$ can be seen 2-simplices, and
so force $g$ to be isomorphic to the formal inverse $f^{-1}$. So 
any word $[f_n]\ldots [f_2][f_1]$ is isomorphic to one not involving any
formal inverses. But now the word $[f_2][f_1]$ is isomorphic to the 
actual composite
$[f_2 f_1]$, and by an easy induction any longer word is likewise
isomorphic to an arrow in \AA.

We now turn to the fullness on 2-cells of $E:\Pi_3N\AA\to\AA$. 
It suffices to show fullness on 2-cells $[f]\to [g]$, where $f$ and
$g$ are 1-cells in \AA.
Any 2-cell $\phi:[f]\to [g]$ in \AA can be regarded as a 2-cell 
$\phi:[f][1]=[f]\to [g]$, 
and so we do indeed have fullness on 2-cells.

The 3-cells of $\Pi_3N\AA$ are once again obtained by adjoining invertible
3-cells 
$$
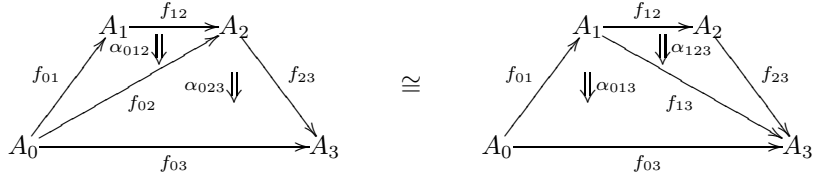
\xy
(0,-8)*{A_0}="A0";
(12,8)*{A_1}="A1";
(28,8)*{A_2}="A2";
(40,-8)*{A_3}="A3";
{\ar^{f_{01}} "A0";"A1"};
{\ar^{f_{12}} "A1";"A2"};
{\ar^{f_{23}} "A2";"A3"};
{\ar_{f_{03}} "A0";"A3"};
{\ar_{f_{02}} "A0";"A2"};
{\ar@{=>}_{\alpha_{023}} (28,2);(28,-2)};
{\ar@{=>}_{\alpha_{012}} (18,7);(18,3)};
\endxy
\qquad \cong \qquad
\xy
(0,-8)*{A_0}="A0";
(12,8)*{A_1}="A1";
(28,8)*{A_2}="A2";
(40,-8)*{A_3}="A3";
{\ar^{f_{01}} "A0";"A1"};
{\ar^{f_{12}} "A1";"A2"};
{\ar^{f_{23}} "A2";"A3"};
{\ar_{f_{03}} "A0";"A3"};
{\ar_{f_{13}} "A1";"A3"};
{\ar@{=>}^{\alpha_{013}} (12,2);(12,-2)};
{\ar@{=>}^{\alpha_{123}} (22,7);(22,3)};
\endxy
$$
% $$\xymatrix{
% & A_1 \ar[r]^{f_{12}} & A_2 \ar[dr]^{f_{23}} && &&
% & A_1 \ar[r]^{f_{12}} \ar[drr]_{f_{13}} \drrtwocell\omit{<-2>\alpha_{123}} & A_2 \ar[dr]^{f_{23}} & \\
% A_0 \ar[ur]^{f_{01}} \ar[urr]_{f_{02}} \ar[rrr]_{f_{03}} \urrtwocell\omit{<-2>\alpha_{012}} & {}\rtwocell\omit{<-2>\alpha_{023}} &{} & A_3 &{}\ar@3[r]^{X}  & {} &
% A_0 \ar[ur]^{f_{01}} \ar[rrr]_{f_{03}} &
% {}\rtwocell\omit{<-2>\alpha_{013}} &{} & A_3 }$$
between the formal composites on either side of the diagram, for 
each 3-cell in \AA between the actual composites. But this time it is
not done freely; rather we introduce an equation for each (commutative) 4-simplex in $N\AA$. (We have not mentioned here the fact that there
will generally already be some 3-cells --- the pseudonaturality 
isomorphisms whose presence is forced whenever there are 2-cells 
whose horizontal composite is not fully determined.)

$\mathbf{\alpha}$

Consider a pair of 2-cells $\bb\alpha,\bb\beta:\bb f\to\bb g$ in $\Pi_3N\AA$.
We wish to show that $E:\Pi_3N\A\to\A$ is fully faithful on 3-cells
from $\bb\alpha$ to $\bb\beta$. We have already seen that $\bb f$ and $\bb g$
can be replaced by isomorphic 1-cells $[f]$ and $[g]$ coming from \AA. 
The fully faithfulness on 3-cells will not be affected by 
whiskering by the isomorphisms $[f]\cong\bb f$ and $[g]\cong\bb g$, so 
we may as well suppose that $\bb\alpha,\bb\beta:[f]\to [g]$. Furthermore,
we can inductively build up isomorphisms $[\alpha]\cong\bb\alpha$ 
and $[\beta]\cong\bb\beta$ in $\Pi_3N\AA$, where $\alpha,\beta:f\to g$ are
now 2-cells in \AA. Once again, it will suffice to prove fully
faithfulness on 3-cells $[\alpha]\to[\beta]$.

Any 3-cell $X:\alpha\to\beta$ in \AA  can be seen as a (fairly degenerate) 
3-simplex
$$
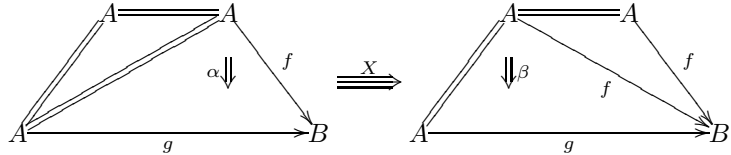
\xymatrix{
{\xy
(0,-8)*{A}="A0";
(12,8)*{A}="A1";
(28,8)*{A}="A2";
(40,-8)*{B}="A3";
{\ar@{=} "A0";"A1"};
{\ar@{=} "A1";"A2"};
{\ar^{f} "A2";"A3"};
{\ar_{g} "A0";"A3"};
{\ar@{=} "A0";"A2"};
{\ar@{=>}_{\alpha} (28,2);(28,-2)};
%{\ar@{=>}_{\alpha_{012}} (15,6);(15,2)};
\endxy} 
\ar@3 [r]^{X} &
%\qquad \cong \qquad
{\xy
(0,-8)*{A}="A0";
(12,8)*{A}="A1";
(28,8)*{A}="A2";
(40,-8)*{B}="A3";
{\ar@{=} "A0";"A1"};
{\ar@{=} "A1";"A2"};
{\ar^{f} "A2";"A3"};
{\ar_{g} "A0";"A3"};
{\ar_{f} "A1";"A3"};
%{\ar@{=>}^{\alpha_{013}} (12,2);(12,-2)};
{\ar@{=>}^{\beta} (12,2);(12,-2)};
\endxy} }
$$
% $$\xymatrix{
% & A \ar[r]^{f} & B \ar@{=}[dr] && &&
% & A \ar[r]^{f} \ar[drr]_{f}  & B \ar@{=}[dr] & \\
% A \ar@{=}[ur] \ar[urr]_{f} \ar[rrr]_{g} &
% {}\rtwocell\omit{<-2>\alpha''} &{} & B &{}\ar@3[r]^{X}  & {} &
% A \ar@{=}[ur] \ar[rrr]_{g} &
% {}\rtwocell\omit{<-2>\beta''} &{} & B }$$
and so can be realized via $E$. This gives fullness on 3-cells of $E$.

Suppose that $\bb X:[\alpha]\to[\beta]$ is any 3-cell in $\Pi_3N\AA$.
Then $\bb X$ can be built up as a formal composite of 3-cells in \AA;
this involves whiskering by 2-cells or 1-cells, and vertical composition
of 3-cells. We can always use the previous arguments to restrict to
the case where
these 3-cells have domain and codomain 2-cells coming from \AA (in 
other words, of the from $[\gamma]$ or $[\delta]$ for 2-cells $\gamma$
or $\delta$ in \AA). 

If $\bb X$ itself has the form $[Y]$ for some
3-cell $Y:\alpha\to\beta$ in \AA, then clearly $X=E\bb X=E[Y]=Y$, and so 
$\bb X=[Y]=[E\bb X]$. We now prove inductively that $\bb X=[E\bb X]$ for
any $\bb X:[\alpha]\to[\beta]$. Any such $\bb X$ has the form 
$[X_n]\ldots[X_1]$; for simplicity we give only the case $\bb X=[X_2][X_1]$;
the general inductive step is essentially the same, but the notation 
becomes more complicated. Write $Y=X_2X_1$, and write $\gamma$ for the
domain of $X_2$ and codomain of $X_1$. Then we are to prove that
$\bb X=[Y]$. 

There is a 4-simplex
$$
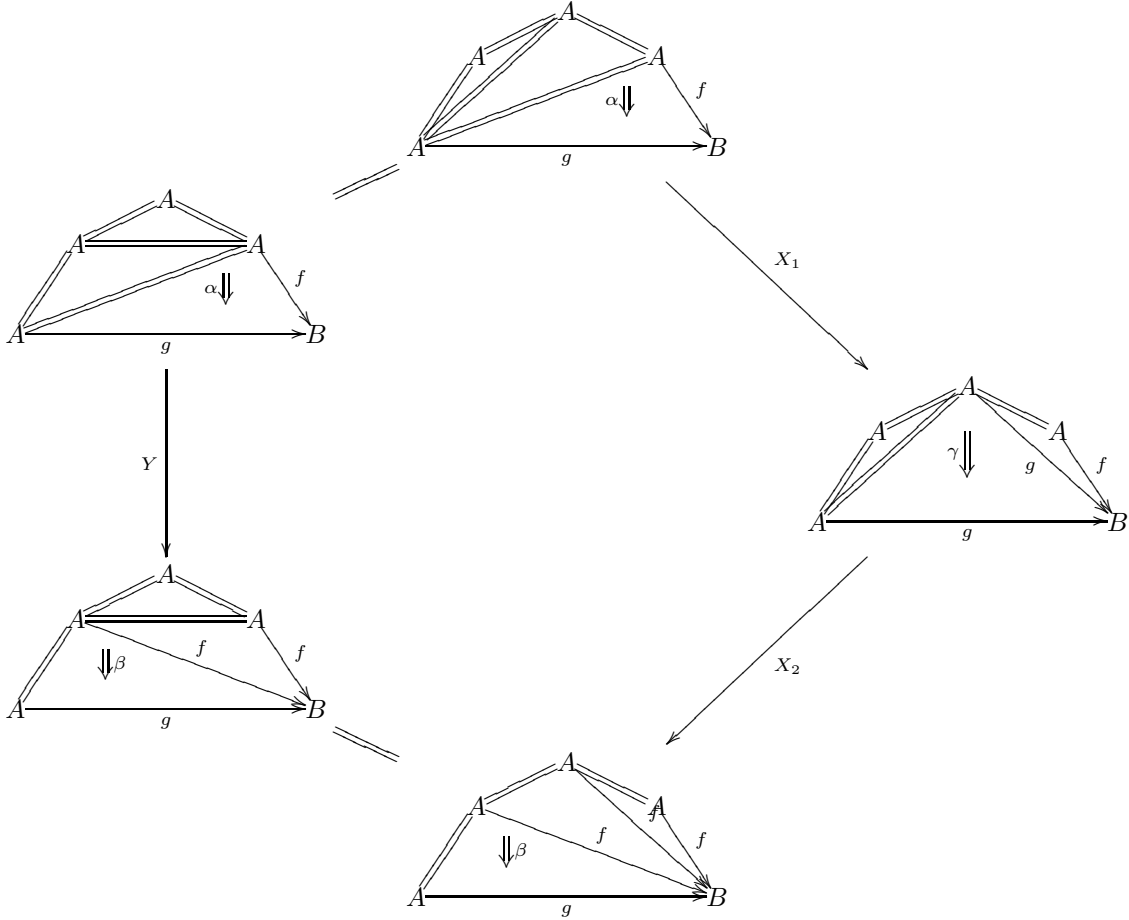
\xymatrix @R0pc { 
& {
\xy
(0,-9)*{A}="A0";
(8,3)*{A}="A1";
(20,9)*{A}="A2";
(32,3)*{A}="A3";
(40,-9)*{B}="A4";
{\ar@{=} "A0";"A1"};
{\ar@{=} "A1";"A2"};
{\ar@{=} "A2";"A3"};
{\ar^{f} "A3";"A4"};
{\ar_{g} "A0";"A4"};
{\ar@{=} "A0";"A2"};
{\ar@{=} "A0";"A3"};
{\ar@{=>}_{\alpha} (28,-1);(28,-5)};
\endxy
}  \ar@{=}[dl] \ar[ddr]^{X_1} \\
{\xy
(0,-9)*{A}="A0";
(8,3)*{A}="A1";
(20,9)*{A}="A2";
(32,3)*{A}="A3";
(40,-9)*{B}="A4";
{\ar@{=} "A0";"A1"};
{\ar@{=} "A1";"A2"};
{\ar@{=} "A2";"A3"};
{\ar^{f} "A3";"A4"};
{\ar_{g} "A0";"A4"};
{\ar@{=} "A1";"A3"};
{\ar@{=} "A0";"A3"};
{\ar@{=>}_{\alpha} (28,-1);(28,-5)};
\endxy
} \ar[dd]_{Y} \\
&&
{\xy
(0,-9)*{A}="A0";
(8,3)*{A}="A1";
(20,9)*{A}="A2";
(32,3)*{A}="A3";
(40,-9)*{B}="A4";
{\ar@{=} "A0";"A1"};
{\ar@{=} "A1";"A2"};
{\ar@{=} "A2";"A3"};
{\ar^{f} "A3";"A4"};
{\ar_{g} "A0";"A4"};
{\ar@{=} "A0";"A2"};
{\ar_{g} "A2";"A4"};
{\ar@{=>}_{\gamma} (20,3);(20,-3)};
\endxy
} \ar[ddl]^{X_2}  \\
{\xy
(0,-9)*{A}="A0";
(8,3)*{A}="A1";
(20,9)*{A}="A2";
(32,3)*{A}="A3";
(40,-9)*{B}="A4";
{\ar@{=} "A0";"A1"};
{\ar@{=} "A1";"A2"};
{\ar@{=} "A2";"A3"};
{\ar^{f} "A3";"A4"};
{\ar_{g} "A0";"A4"};
{\ar@{=} "A1";"A3"};
{\ar^{f} "A1";"A4"};
{\ar@{=>}^{\beta} (12,-1);(12,-5)};
\endxy
} \ar@{=}[dr] \\
&
{\xy
(0,-9)*{A}="A0";
(8,3)*{A}="A1";
(20,9)*{A}="A2";
(32,3)*{A}="A3";
(40,-9)*{B}="A4";
{\ar@{=} "A0";"A1"};
{\ar@{=} "A1";"A2"};
{\ar@{=} "A2";"A3"};
{\ar^{f} "A3";"A4"};
{\ar_{g} "A0";"A4"};
{\ar^{f} "A2";"A4"};
{\ar^{f} "A1";"A4"};
{\ar@{=>}^{\beta} (12,-1);(12,-5)};
\endxy
}
}$$
\noindent
in $N\AA$, and so  in $\Pi_3N\AA$ there is a relation $[X_2][X_1]=[Y]$,
which is just to say $\bb X=[E\bb X]$. 

This now proves the faithfulness on 3-cells.
\endproof

\begin{corollary}\label{cor:transport}
The model structure on \graygpd is obtained by transporting along the nerve
functor $N:\graygpd\to\SSet$ the model structure on \SSet.
\end{corollary}

\proof
We know that $N$ preserves weak equivalences, fibrations, and trivial
fibrations; in other words that if a \Gray-functor $F:\AA\to\BB$ is a weak
equivalence, fibration, or trivial fibration, then so is $NF$. It remains to 
show the converse of each of these three statements; in fact, it will suffice
to do two of the three.

Suppose first that $NF$ is a weak equivalence in \SSet. Then $\Pi_3 NF$
is a weak equivalence in \graygpd; but $E:\Pi_3 N\AA\to\AA$ and 
$E:\Pi_3 N\BB\to\BB$ are weak equivalences, then so is $F$.

Now suppose that $NF$ is a trivial fibration in \SSet. Then certainly
it has the right lifting property with respect to the inclusion 
$\partial\Delta_0\to\Delta_0$; and so by adjointness $F$ has the 
right lifting property with respect to $\Pi_3\partial\Delta_0\to\Pi_3\Delta_0$;
but this is just the unique map from the initial \Gray-groupoid to the terminal
one, so the right lifting property says that $F$ is surjective on objects.

Similarly, $NF$ has the right lifting property with respect to 
$\partial\Delta_1\to\Delta_1$, and so $F$ is full on 1-cells.

Fullness on 2-cells and 3-cells is slightly more complicated, since 
\Gray-categories have a globular rather than a simplicial structure. We can 
build a simplicial model ``$\two_2$'' of a parallel pair
$$
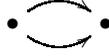
\xymatrix{\bullet \rtwocell{\omit} & \bullet}$$
by taking the coproduct $\Delta_1+\Delta_1$ and identifying the 
boundaries. Similarly we can build a simplicial model ``$\two_\two$'' 
of a (globular) 2-cell
$$
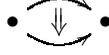
\xymatrix{\bullet \rtwocell & \bullet}$$
by collapsing one of the faces of a 2-simplex. The inclusion $\two_2\to\two_\two$
is a monomorphism, so its image under $\Pi_3$ is a cofibration; it follows that 
$F$ is full on 2-cells. 

Fullness on 3-cells is once again proved by a constructing simplicial model of
a parallel pair of 2-cells, and a simplicial model of a 3-cell, and observing that the
inclusion is a monomorphism.

Finally, faithfulness on 3-cells can be interpreted as ``fullness on 4-cells'':
think of a Gray-category as some sort of 4-dimensional category in which there 
are no non-identity 4-cells; then equality of 3-cells becomes the existence of 
a (necessarily trivial) 4-cell between them. Let $f,g:A\to B$ be 1-cells in \AA, with 
2-cells $\alpha,\beta:f\to g$
and 3-cells $M,M':\alpha\to\beta$. We need to show that $FM=FM'$ implies that 
$M=M'$; in other words that any 4-cell $FM\to FM'$ is in the image under $F$ 
of a 4-cell $M\to M'$. This can once again be expressed simplicially. 
\endproof

\section{Localizations of \SSet}
\label{sect:localization}

In this section we describe the relationship between $n$-groupoids and 
homotopy $n$-types (for small $n$) in terms of Bousfield localizations
\cite{Hirschhorn-book} of 
\SSet. For each $n$, let $g_n$ be the inclusion $\partial\Delta_{n+2}\to\Delta_{n+2}$.
The localizations we consider will be with respect to these maps $g_n$.
(Once could also use the image of the inclusion $S^{n+1}\to D^{n+2}$ in \Top
under the singular functor $\Top\to\SSet$ as $g_n$.)

We start by recalling the well-understood situation for $n=1$. The 
nerve functor $N:\gpd\to\SSet$ for groupoids has a left adjoint 
$\Pi_1:\SSet\to\gpd$ sending a simplicial set to its fundamental groupoid 
\cite{GabrielZisman}. The nerve functor preserves fibrations and trivial 
fibrations, so is a right Quillen functor. The left adjoint $\Pi_1$ therefore
preserves all cofibrations and trivial cofibrations. But since all objects of \SSet are 
cofibrant, and all objects of \gpd are fibrant, in fact both functors preserve
weak equivalences. The fundamental groupoid $\Pi_1(X)$ of a simplicial set
agrees with the ordinary fundamental groupoid $\Pi_1|X|$ of the geometric
realization of $X$.

Now the nerve functor is fully faithful, so the counit $\Pi_1 N\to 1$ is 
invertible, and so $\Pi_1$ can be thought of as a reflection onto a full
subcategory. The unit is not invertible, and $\Pi_1$ does destroy information:
in particular, it destroys all homotopical information in dimension greater than 1.
More precisely, for a morphism $f:X\to Y$ of simplicial sets, the induced
functor $\Pi_1(f):\Pi_1(X)\to\Pi_1(Y)$ is an equivalence if and only if 
$f$ induces a bijection $\pi_0 f:\pi_0 X\to\pi_0 Y$ between the sets of 
connected components, and an isomorphism $\pi_1 f:\pi_1(X,x)\to\pi_1(Y,fx)$
of fundamental groups for all choices of basepoint $x\in X$. 

Recall that $g_1$ is the inclusion $\partial\Delta_3\to\Delta_3$.
%  image, under the singular functor $\Top\to\SSet$, of the 
% inclusion $S^2\to D^3$.
Let $P_1\SSet$ be the Bousfield localization 
of \SSet with respect to the map $g_1$. The weak equivalences are
precisely
the morphisms $f:X\to Y$ of simplicial sets for which
$\Pi_1(f):\Pi_1(X)\to\Pi_1(Y)$ is an equivalence, and the cofibrations
are the usual cofibrations of simplicial sets (the
monomorphisms). Thus the Quillen 
adjunction between \gpd and \SSet passes
to a Quillen adjunction between \gpd and $P_1\SSet$, which is now a
Quillen  equivalence.

Next we turn to the case $n=2$. Once again there is a nerve functor
$N:\twogpd\to\SSet$. This nerve of a 2-groupoid is a special case of
Street's nerve of a bicategory. As observed in
\cite{MoerdijkSvensson}, this nerve functor $N:\twogpd\to\SSet$ is the 
right adjoint part of a Quillen adjunction. We shall write $\Pi_2$ for
the left adjoint, called $W$ in \cite{MoerdijkSvensson}. Once again
$N$ and $\Pi_2$ both preserve weak equivalences; a simplicial map 
$f:X\to Y$ induces a weak equivalence of 2-groupoids if and only if 
it induces a bijection $\pi_0(f):\pi_0(X)\to\pi_0(Y)$ of sets, and 
isomorphisms $\pi_1(f):\pi_1(X,x)\to\pi_1(Y,fx)$ and
$\pi_2(f):\pi_2(X,x)\to\pi_2(Y,y)$ of groups for all basepoints $x\in
X$.  

Let $P_2\SSet$ be the Bousfield localization 
of \SSet with respect to the map $g_2$. This has the same cofibrations 
as \SSet, and $f:X\to Y$ is a weak equivalence in $P_2\SSet$ if and
only if $\Pi_2(f):\Pi_2(X)\to\Pi_2(Y)$ is one in \twogpd. The Quillen
adjunction $\Pi_2\dashv N:\twogpd\to\SSet$ passes
to a Quillen equivalence between \twogpd and $P_2\SSet$.

Finally, we turn to the case $n=3$. We saw in the previous section
that the nerve functor $N:\graygpd\to\SSet$ of \cite{Berger-3type} is the
right adjoint part of a Quillen adjunction; we write $\Pi_3$ for the 
left adjoint part. Once again
$N$ and $\Pi_3$ both preserve weak equivalences; a simplicial map 
$f:X\to Y$ induces a weak equivalence of Gray-groupoids if and only if 
it induces a bijection $\pi_0(f):\pi_0(X)\to\pi_0(Y)$ of sets, and 
isomorphisms $\pi_n(f):\pi_n(X,x)\to\pi_n(Y,fx)$ of groups for all
basepoints $x\in X$, and all $n\in\{1,2,3\}$.

Let $P_3\SSet$ be the Bousfield localization 
of \SSet with respect to the map $g_3$. This has the same cofibrations 
as \SSet, and $f:X\to Y$ is a weak equivalence in $P_3\SSet$ if and
only if $\Pi_3(f):\Pi_3(X)\to\Pi_3(Y)$ is one in \graygpd. The Quillen
adjunction $\Pi_3\dashv N:\graygpd\to\SSet$ passes
to a Quillen equivalence between \graygpd and $P_3\SSet$.

\section{Tricategories}\label{sect:tricategories}

We conjecture that the model structure on \graycat can be extended to the context of tricategories. To do this, one should use a ``fully algebraic''
definition of tricategories, such as that of Gurski
\cite{Gurski-thesis}, 
so that the category \tricat
of tricategories and strict morphisms of tricategories is locally presentable,
and contains \graycat as a full reflective subcategory.

The definitions of weak equivalence and fibration still make perfectly
good sense in the context of tricategories; indeed the notion of 
triequivalence (our weak equivalences) was first defined in that context
\cite{GPS}. Choosing these weak equivalences and fibrations determines 
a model structure uniquely if one exists. We conjecture that it does.
The inclusion of \graycat into \tricat preserves fibrations and weak
equivalences, and so would be the right adjoint part of a Quillen adjunction.

We conjecture further that this Quillen adjunction is in fact a Quillen
equivalence. Since the inclusion $I:\graycat\to\tricat$ is fully faithful,
the counit $LI\to 1$ of the adjunction is invertible. Since all objects
of \graycat are fibrant, the unit of the derived adjunction will also be
invertible; so it would remain only to show that the counit of the derived
adjunction is invertible. This will be the case provided that the 
counit $\TT\to IL\TT$ is a weak equivalence for every cofibrant tricategory
\TT, and Michael Makkai has announced \cite{Makkai-talk} that he has
proved this to be the case.

On the other hand, we could consider strict 3-categories.
The category \threecat of 3-categories and 3-functors is a full reflective subcategory of \graycat. The model structure on \graycat lifts across the
adjunction to give a model structure on \threecat for which a 3-functor
is a weak equivalence or a fibration if and only if the corresponding 
\Gray-functor is one. One simply applies the reflection into \graycat
to each of the generating cofibrations and generating trivial cofibrations
(in fact this is only needed for one generating trivial cofibration; the
other generating trivial cofibrations and the generating cofibrations 
are already in \threecat) and then observes that for a 3-category \BB,
the path object \PB in \graycat constructed above is in fact a 3-category. 

The inclusion of \threecat in \graycat together with its left adjoint 
are then a Quillen adjunction, by construction, but they will not be a
Quillen equivalence, essentially because not every 3-category is
triequivalent to a 3-category \cite{GPS}.

\section{Computads and sesquicategories}
\label{sect:computads}

We have given explicit descriptions of the fibrations and weak 
equivalences, but so far the cofibrations have been defined only via a left lifting property. 
In the following section we investigate the cofibrations of \graycat, and in particular, the cofibrant objects. First, however, we develop some necessary
material on sesquicategories \cite{Street:categorical-structures}.

First, recall that the functor $U:\twocat\to\cat$ which forgets the
2-cells of a 2-category has both a left adjoint $D$, which adds
identity 2-cells to a category, and a right adjoint $C$, which adds a
single (invertible) 2-cell between each parallel pair of 1-cells. Both $D$ and $C$ are fully faithful, so the 
unit $1\to UD$ of the adjunction $D\dashv U$ and the counit $UC\to 1$ of 
the adjunction $U\dashv C$ are invertible.

We shall regard the category \twocat as monoidal via the Gray tensor product, and \cat as monoidal via the ``funny'' tensor 
product $\square$. Recall, for example from \cite{Street:categorical-structures}, 
 that $\square$ is part of a monoidal closed structure
in which the internal hom $[A,B]$ has as objects the functors from $A$ to 
$B$ and as morphisms from $f$ to $g$, the set of all $\ob A$-indexed families
of morphisms $fa\to ga$, with no naturality condition imposed. 

This determines the funny tensor product; it can be described more
explicitly as follows: its objects are pairs $(a,b)$ of objects of $A$
and $B$, while its morphisms are freely generated by morphisms of the form $(\alpha,1_b)$ and
$(1_a,\beta)$, where $\alpha:a\to a'$ is a morphism in $A$ and $\beta:b\to b'$
a morphism in $B$; subject to relations which state that $(a,-)$ and $(-,b)$ are functors. 
Another way to say this is $A\square B=U(CA\ox CB)=U(DA\ox DB)$; and in fact 
both $U$ and $C$ are strong monoidal functors, while $D$ is only opmonoidal.

Now a category enriched in \twocat (with respect to the Gray tensor product)
is of course a Gray-category, while a category enriched in \cat (with respect
to the funny tensor product) is called a {\em sesquicategory}. Since the
adjunction $U\dashv C$ is monoidal, we obtain an adjunction $U_*\dashv C_*$
between \graycat and \sesquicat, where $U_*$ discards the 3-cells, and $C_*$ 
adjoins a single (invertible) 3-cell between each parallel pair of 2-cells. Since $D$
is not monoidal it does not induce a 2-functor $D_*$, but there is 
nonetheless a left adjoint $L$ to $U_*$, which adjoins only identity 2-cells
and pseudonaturality isomorphisms. Note that the counits $LU_*\A\to\A$ 
of the adjunction $L\dashv U_*$ are all bijective on objects, 1-cells, and
2-cells; they only change the 3-cells.

We shall also need a notion of {\em free} sesquicategory on a
computad. Recall, for example from \cite{Street:categorical-structures},
that a computad
consists of following structure: first of all a directed graph
$d,c:G_1\to G_0$, on which we form the free category $C$; then a
second directed graph $d,c:G_2\to C_1$, whose vertices are the 1-cells
of $C$; finally we require that the ``globular relations'' $dd=dc$
and $cd=cc$ hold, so that elements of $G_2$ (called 2-cells) can be
written as the diagram below on the left where $f,g:x\to y$ are 1-cells in $C$
or in expanded from as in the diagram on the right
$$
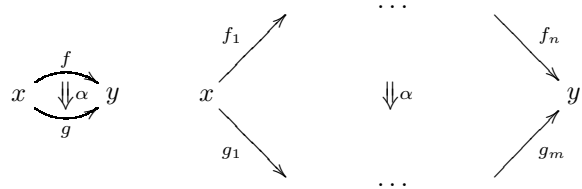
\xymatrix{
& &&& {}  \ldots  & \relax\ar[dr]^{f_{n}} \\
x \rtwocell^f_g{\alpha} & y & x \ar[ur]^{f_1} \ar[dr]_{g_1} 
& \relax\rrtwocell\omit{\alpha} && \relax & y \\
& &&& {}  \ldots & \relax\ar[ur]_{g_{m}}  }$$
in which the ellipses denote an unspecified number of further
composites (the picture on the right appears to suggest that $n\ge 2$,
but in fact $n$ could be $0$ or $1$).
A morphism from such a computad $G$ to a computad $H$ consists of a graph morphism $f$ from 
$d,c:G_1\to G_0$ to $d,c:H_1\to H_0$, consisting of $f_0:G_0\to H_0$
and $f_1:G_1\to H_1$ compatible with the source and target maps; equipped
with a map $f_2:G_2\to H_2$ satisfying the evident compatibility condition.

These form a category \cmptd which, as observed by Schanuel, is in
fact a presheaf category. There is a forgetful functor $V:\sesquicat\to\cmptd$
which sends a sesquicategory \A to the computad whose 0-cells and 1-cells
are the 0-cells and 1-cells of \A, and whose 2-cells between the paths
$f_n\ldots f_2f_1\to g_m\ldots g_2g_1$ are the 2-cells between their 
composites in \A. 

This forgetful functor $V$ has a left adjoint $H\dashv V$, which can
be constructed as follows. An object of $HG$ is an element of $G_0$.
A 1-cell of $HG$ is a path in the 1-cells of $G$. A ``basic 2-cell'',
is an expression of the form $\ell\alpha r$, where $\ell:w\to x$ and $r:y\to z$
are 1-cells in $HG$, and $\alpha:f\to g:x\to y$ is an element of $G_2$. 
Such a basic 2-cell has a source 1-cell $\ell fr$ and target 1-cell
$\ell gr$. A 2-cell of $HG$ is a path in the basic 2-cells.

We shall call a sesquicategory {\em free}, if it is (isomorphic to) one
of the form $HG$ for a computad $G$. 
In \cite[Lemma~4.7]{qm2cat} it 
was proved that any retract of a free category is free. We use this 
fact to prove an analogue for sesquicategories.

\begin{proposition}\label{prop:retracts}
A retract of a free sesquicategory is free.
\end{proposition}

\proof
Let $G$ be a computad, $HG$ the free sesquicategory on $G$, and
let $I:\A\to HG$ be a sesquifunctor with a retract $R:HG\to\A$.
Certainly the underlying category of $HG$ is free, and the underlying
category of \A is a retract of it, so the underlying category of \A
is also free. 

As well as the underlying category of a sesquicategory, there is 
also the category obtained by discarding the objects: it has arrows
of the sesquicategory as objects, 2-cells of the sesquicategory 
as arrows, and vertical composition in the sesquicategory as composition.
This clearly induces a functor $W:\sesquicat\to\cat$. If \A is a 
retract of $HG$ then $W\A$ is a retract of $WHG$; but $WHG$ is the 
free category on all whiskerings of 2-cells in $HG$, so $W\A$ is also 
free. 

It remains to show that the generating arrows of $W\A$ are freely 
generated under whiskering in \A by some subset. Say that a 2-cell
$\alpha:f\to g$ in \A is indecomposable if it is indecomposable as
an arrow of $W\A$: these indecomposable 2-cells are the generating
arrows of $W\A$. Say that such an indecomposable 2-cell is
 {\em totally indecomposable} if
moreover it cannot be written as a whiskering $\alpha=g\beta f$ 
unless $f$ and $g$ are both identities. We need to show that every
indecomposable 2-cell can be written uniquely as a whiskering
of a totally indecomposable 2-cell.

Let $\alpha:f\to g$ be any 2-cell in \A, and consider its image
$I\alpha:If\to Ig$ in $HG$. We can write $I\alpha$ uniquely as 
a composite $\beta_n\ldots\beta_1$ of indecomposables; and each
of these can in turn be written uniquely as a whiskering
$\beta_i=r_i\gamma_i \ell_i$ with $\gamma_i$ totally indecomposable.
Since the underlying category of $HG$ is free, each $r_i$ and each
$\ell_i$ has a well-defined length; write $\pi(\alpha)$ for the sum
of all these lengths.

If there is an indecomposable 2-cell $\alpha':f'\to g':x'\to y'$ 
which is not the whiskering of any totally indecomposable 2-cell,
then there is one with $\pi(\alpha')$ minimal.
Then certainly $\alpha'$ is not totally indecomposable, so we can
write it as a whiskering $\alpha'=r\alpha\ell$, where $\alpha:f\to g:x\to y$,
and $r$ and $\ell$ are not both identities.
Decompose $I\alpha$ as above: $I\alpha=\beta_n\ldots\beta_1$ where
$\beta_i=r_i\gamma_i\ell_i$. But now the decomposition of $I\alpha'$ 
is $(r\beta_n\ell)\ldots(r\beta_1\ell)$, with $r\beta_i\ell=rr_i\beta_i\ell_i\ell$; and so $\pi(\alpha)<\pi(\alpha')$
contradicting the minimality of $\pi(\alpha')$.

This proves that every indecomposable 2-cell is the whiskering of
some totally indecomposable 2-cell. It remains to show the uniqueness.
Suppose then that $\alpha:f\to g:x\to y$ and $\alpha':f'\to g':x'\to y'$
are totally indecomposable 2-cells in \A, and that 
$r\alpha\ell=r'\alpha'\ell'$ for some 1-cells $r:y\to z$, $r':y'\to z$,
$\ell:w\to x$, and $\ell':w\to x'$.

Write $I\alpha=\beta_n\ldots\beta_1$ and $\beta_i=r_i\gamma_i\ell_i$
where each $\gamma_i$ is totally indecomposable in $HG$; and similarly
$I\alpha'=\beta'_{n'}\ldots\beta'_1$ and $\beta'_i=r'_i\gamma'_i\ell'_i$.
Now 
$Ir.I\alpha.I\ell=I(r)r_n\gamma_n\ell_nI(\ell)\ldots I(r)r_1\gamma_1\ell_1I(\ell)$ and 
$Ir'.I\alpha'.I\ell'=I(r')r'_{n'}\gamma'_{n'}\ell'_{n'}I(\ell')\ldots
I(r')r'_1\gamma'_1\ell'_1I(\ell')$ are both decompositions into
totally indecomposables of the same 2-cell in $HG$, and so must agree. 
Thus $n=n'$ and $\gamma_i=\gamma'_i$ for each $i$, while also 
$I(r)r_i=I(r')r'_i$ and $\ell_iI(\ell)=\ell'_iI(\ell')$ for each $i$.

Since $I(r)r_i=I(r')r'_i$, either there is a 1-cell $s$ such 
that $I(r)=I(r')s$ and $sr_i=r'_i$, or there is a 1-cell
$s$ such that $I(r')=I(r)s$ and $sr'_i=r_i$. Without loss of generality
we take the former.
% Similarly we can use the fact that $\ell_i I(\ell)=\ell'_i I(\ell')$.
% Once again there are two possibilities; we suppose this time 
% $I(\ell)=mI(\ell')$ and $\ell_i m=\ell'_i$ (the other possibility 
% $I(\ell')=mI(\ell)$ and $\ell'_i m=\ell_i$ is similar). 
Applying $R$ we get $r'=rR(s)$, and now
\begin{align*}
\alpha' &= RI\alpha' \\
       &= R(\beta'_n)\ldots R(\beta'_1) \\
       &= R(r'_n\gamma_n\ell'_n)\ldots R(r'_1\gamma'_1\ell'_1) \\
       &= R(sr_n\gamma_n\ell'_n)\ldots R(sr_1\gamma_1\ell'_1) \\
       &= Rs.\bigl(R(r_n\gamma_n\ell'_n)\ldots R(r_1\gamma_1\ell'_1)\bigr) 
\end{align*}
but $\alpha'$ was assumed totally indecomposable, and so $Rs$ must
be an identity, $r=r'$.

A similar argument shows that $\ell=\ell'$. Now
$$Ir.I\alpha.I\ell=Ir'.I\alpha'.I\ell'$$
in $HG$, with $r=r'$ and $\ell=\ell'$, and so $I\alpha=I\alpha'$
since $HG$ is free, and finally $\alpha=\alpha'$.
\endproof

\section{Cofibrations}
\label{sect:cofibrations}

The first step of our analysis of cofibrations in \graycat is to reduce the
problem to one about sesquicategories. To this end, we define a 
sesquifunctor $F:\A\to\B$ to be a {\em surjection} if it is surjective on 
objects and on each hom $F:\A(A,B)\to\B(FA,FB)$ it is a surjective equivalence
(a trivial fibration in \cat). We also define a sesquicategory \C to be
{\em projective} if the function $\sesquicat(\C,F):\sesquicat(\C,\A)\to\sesquicat(\C,\B)$ is surjective 
for all surjective sesquifunctor $F:\A\to\B$.

The next result is a direct analogue
of \cite[Lemma~4.1]{qm2cat}.

\begin{proposition}
A \Gray-functor $F:\A\to\B$ is a cofibration if and only if its underlying sesquifunctor $U_*F:U_*\A\to U_*\B$  has the left lifting property with 
respect to surjective sesquifunctors. In particular, a \Gray-category \A
is cofibrant if and only its underlying sesquicategory $U_*\A$ is 
projective (with respect to the surjective sesquifunctors).
\end{proposition}

\proof
If $P:\C\to\D$ is a surjective sesquifunctor then $C_*P:C_*\C\to C_*\D$ is a 
trivial fibration of \Gray-categories. If $F$ is a cofibration then it will
have the left lifting property with respect to $C_*P$, and so $U_*F$ will
have the left lifting property with respect to $P$. 

This proves one direction; for the other, suppose that $U_*F$ does have the
left lifting property with respect to surjective sesquifunctors. This time,
let $P:\C\to\D$ be a trivial fibration of \Gray-categories. Then $U_* P$
is clearly surjective, and so $U_* F$ has the left lifting property with 
respect to $U_* P$; but this now means that $L U_* F$ has the left lifting
property with respect to $P$. If now $X$ and $Y$ are \Gray-functors 
satisfying $PX=YF$ we have the following diagram of \Gray-functors
$$\xymatrix{
LU_*\A \ar[r]^{I} \ar[d]_{LU_* F} & \A \ar[r]^{X} \ar[d]_{F} & \C \ar[d]^{P} \\
LU_*\B \ar[r]_{J} & \B \ar[r]_{Y} & \D }$$
in which $I$ and $J$ denote the canonical inclusions (components of the 
counit of $L\dashv U_*$). Since $LU_*F$ has the left lifting property with
respect to $P$, there is a \Gray-functor $Z:LU_*\B\to\C$ whose composite
with $LU_*F$ is $XI$ and whose composite with $P$ is $YJ$. Now $J$ is 
bijective on objects, 1-cells, and 2-cells; while $P$ is locally locally 
fully faithful (bijective on 3-cells with given domains and codomains), 
and so there is a unique induced $W:\B\to\C$ with $PW=Y$ and $WJ=Z$. 
Finally $PWF=YF=PX$ and $WFI=WJ.LU_*F=Z.LU_*F=XI$, while $I$ is 
bijective on objects, 1-cells, and 2-cells, and $P$ is locally locally 
fully faithful, and so $WF=X$. Thus $X$ provides the desired fill-in
satisfying $WF=X$ and $PW=Y$. \endproof

For a sesquicategory \A, the 
counit $HV\A\to\A$ is bijective on objects, and surjective on 
1-cells. Given 1-cells $f,g:x\to y$ (paths in the 1-cells of \A),
a 2-cell in \A between their composites is included in $HV\A$ as
one of the generating 2-cells. Thus the counit $E:HV\A\to\A$ is 
also full on 2-cells and we have:

\begin{proposition}
For each sesquicategory \A, the counit map $HV\A\to\A$ is a 
surjective sesquifunctor. \endproof
\end{proposition}

We now find ourselves in the typical situation where projectives
are the retracts of frees.

\begin{proposition}
For a sesquicategory \A the following are equivalent:
\begin{enumerate}[(i)]
\item \A is projective;
\item there is a sesquifunctor $J:\A\to HV\A$ with $EJ=1$;
\item \A is a retract of a free sesquicategory $HG$ on some computad $G$.
\item \A is free;
\end{enumerate}
\end{proposition}

\proof
The implication $(i)\Rightarrow(ii)$ follows immediately from the fact
that $E$ is a surjection of sesquicategories; the implication
$(ii)\Rightarrow(iii)$ is trivial, and the implication
$(iii)\Rightarrow(iv)$ is Proposition~\ref{prop:retracts}. 
For the implication
$(iv)\Rightarrow(i)$ 
% first observe that retracts of projective sesquicategories are projective,
% and so it suffices to show that each free sesquicategory $HG$ if
% projective. 
we must show that any free sesquicategory $HG$ is projective.
But to say that $HG$ is projective with respect to a surjective
sesquifunctor $P:\A\to\B$ is equivalent to saying that $G$ is 
projective with respect to the underlying computad morphism 
$VP:V\A\to V\B$. But for a surjective sesquifunctor $P:\A\to\B$,
clearly $VP:V\A\to V\B$ has a section, and so any object $G$ will
be projective with respect to $VP$. 
\endproof

\begin{corollary}
A \Gray-category is cofibrant if and only if its underlying sesquicategory
is free on a computad. \endproof
\end{corollary}

We can now give a very explicit description of a cofibrant replacement
functor. For a \Gray-category \AA, first forget the 3-cells, to 
obtain a sesquicategory $U_*\AA$, then form the free sesquicategory $HVU_*\AA$ 
on its underlying computad, and then the free \Gray-category 
$LHVU_*\AA$ on that. This has a canonical \Gray-functor 
$E':LHVU_*\AA\to\AA$ to the original \Gray-category; in fact this is the 
counit at \AA of an adjunction between \Gray-categories and computads.
This \Gray-functor $E'$ is bijective on objects, surjective on 1-cells,
and full on 2-cells. Factorize it as 
$$\xymatrix{
LHVU_*\AA \ar[r]^{i} & Q\AA \ar[r]^{q} & \AA }$$
where $i$ is bijective on objects, 1-cells, and 2-cells, and $q$ is
fully faithful on 3-cells. Clearly this defines a functor 
$Q:\graycat\to\graycat$ and a natural transformation $q:Q\to 1$.

Since $i$ is bijective on objects, 1-cells,
and 2-cells, and the underlying sesquicategory of $LHVU_*\AA$ is 
just the free sesquicategory $HVU_*\AA$; the underlying sesquicategory 
of $Q\AA$ is also just $HVU_*\AA$, so $Q\AA$ is indeed a cofibrant
\Gray-category. Furthermore, since $E'$ is bijective on objects, surjective
on 1-cells, and full on 2-cells, so is $q$; since $q$ is also fully faithful
on 3-cells, it is in fact a trivial fibration of \Gray-categories.

But in fact we can do a bit better: $Q$ is actually a comonad, so we can 
obtain a category of ``weak'' morphisms of \Gray-categories by taking 
the Kleisli category. This turns out to be a special case of a general
construction due to Garner \cite{Garner:homomorphisms}.

We already have the counit $q:Q\to 1$; next we build the comultiplication
$d:Q\to Q^2$. To do this, write $W:\graycat\to\cmptd$ for the composite
of $U_*:\graycat\to\sesquicat$ and $V:\sesquicat\to\cmptd$, and $K$
for the left adjoint $HV$; with unit $\eta:1\to VK$ and counit
$\epsilon:KV\to 1$. Observe that the exterior of 
$$\xymatrix{
& FU\AA \ar[dr]^{1} \ar[dl]_{F\eta U\AA} \\
FUFU\AA \ar[rr]_{\epsilon FU\AA} \ar[d]_{FUi\AA} && 
FU\AA \ar[d]^{i\AA} \\
FUQ\AA \ar[rr]^{\epsilon Q\AA} \ar[dr]_{iQ\AA} && Q\AA \\
& QQ\AA \ar[ur]_{qQ\AA} }$$
commutes, and that $i\AA$ is bijective on objects, 1-cells, and 2-cells,
while $qQ\AA$ is fully faithful on 3-cells. Thus there is a unique
\Gray-functor $d:Q\AA\to Q^2\AA$ making the diagram
$$\xymatrix{
FU\AA \ar[r]^{i\AA} \ar[d]_{F\eta U\AA} & Q\AA \ar[ddd]^{1} \ar[dddl]_{d\AA} \\
FUFU\AA \ar[d]_{FUi\AA} \\
FUQ\AA \ar[d]_{iQ\AA} \\
Q^2\AA \ar[r]_{qQ\AA} & Q\AA }$$
commute. These $d\AA$ are the components of a natural transformation
$d:Q\to Q^2$; naturality of $d$ is inherited from that of $\eta$, $\epsilon$,
$i$, and $q$. One of the counit laws $qQ.d=1$ holds by definition of $d$.
We prove the other law $Qq.d=1$, by checking that the composite of each
side with $i$ agrees (so that $Qq.d$ and $1$ agree on objects, 1-cells, and
2-cells), and that the composite of each side with $q$ agrees (so that the
two sides agree on 3-cells). The calculations are:
\begin{align*}
q.Qq.d &= q.qQ.d \\
       &= q.1 \\
Qq.d.i &= Qq.iQ.FUi.F\eta U \\
       &= i.Fuq.FUi.F\eta U \\
       &= i.FU\epsilon.F\eta U \\
       &= i \\
       &= 1.i
\end{align*}
and so $Qq.d=1$. We prove the coassociative law $Qd.d=dQ.d$ using the 
same technique:
\begin{align*}
qQ^2.Qd.d &= d.qQ.d \\
          &= d \\
          &= qQ^2.dQ.d \\
Qd.d.i &= Qd.iQ.FUi.F\eta U \\
       &= iQ^2.FUd.FUi.F\eta U \\
       &= iQ^2.FUiQ.FUFUi.FUF\eta U.F\eta U \\
       &= iQ^2.FUiQ.F\eta UQ.FUi.F\eta U \\
       &= dQ.iQ.FUi.F\eta U \\
       &= dQ.d.i
\end{align*}
and so $Qd.d=dQ.d$.

We have now proved:

\begin{theorem}
There is a comonad $(Q,q,d)$ on \graycat for which $q\AA:Q\AA\to\AA$
is a cofibrant replacement of \AA for every \Gray-category \AA.
\end{theorem}

As mentioned above, this comonad can be obtained using the results of 
\cite{Garner:homomorphisms}, which provide a general technique for
defining weak morphisms of higher categories via Kleisli categories
for comonads.

\bibliographystyle{plain}
%\bibliography{my}

\end{document}